\documentclass{elsarticle}
\pdfoutput=1

\newcommand{\eqpunct}[1]{\text{#1}}
\usepackage{algorithm}
\usepackage{algorithmicx}
\usepackage[noend]{algpseudocode}
\usepackage{pgfplots}
\usepackage{subcaption}
\usepackage{booktabs}
\usepackage{dcolumn}
\usepackage{hyperref}
\usepackage{amsmath}
\usepackage{natbib}

\graphicspath{{./figures_bw/}{./figures/}}
\pgfplotsset{compat=newest}

\def\Ncal{\mathcal{N}}

\def\Ical{\mathcal{I}}

\def\Xcal{\mathcal{X}}

\newcommand{\argmax}{\operatornamewithlimits{arg\,max}}

\begin{document}

\begin{frontmatter}

\title{A comparison of Monte Carlo tree search and mathematical optimization for large scale dynamic resource allocation \tnoteref{t1}}

\tnotetext[t1]{This work is sponsored by the Assistant Secretary of Defense for Research and Engineering, ASD(R\&E), under Air Force Contract \#FA8721-05-C-0002. Opinions, interpretations, conclusions, and recommendations are those of the authors and are not necessarily endorsed by the United States Government.}

\author[MITSloanORC]{Dimitris Bertsimas}
\ead{dbertsim@mit.edu}

\author[LL]{J. Daniel Griffith}
\ead{dan.griffith@ll.mit.edu}

\author[MITORC]{Vishal Gupta}
\ead{vgupta1@mit.edu}

\author[StanfordAA]{Mykel J. Kochenderfer}
\ead{mykel@stanford.edu}

\author[MITORC]{Velibor V. Mi\v{s}i\'{c}}
\ead{vvmisic@mit.edu}

\author[LL]{Robert Moss}
\ead{robert.moss@ll.mit.edu}

\address[MITSloanORC]{Sloan School of Management and Operations Research Center, Massachusetts Institute of Technology; 77 Massachusetts Avenue, Cambridge MA 02139}
\address[LL]{Lincoln Laboratory, Massachusetts Institute of Technology; 244 Wood Street, Lexington MA 02420}
\address[MITORC]{Operations Research Center, Massachusetts Institute of Technology; 77 Massachusetts Avenue, Cambridge MA 02139}
\address[StanfordAA]{Department of Aeronautics and Astronautics, Stanford University; 496 Lomita Mall, Stanford CA 94305}

\begin{abstract}
Dynamic resource allocation (DRA) problems are an important class of dynamic stochastic optimization problems that arise in a variety of important real-world applications. DRA problems are notoriously difficult to solve to optimality since they frequently combine stochastic elements with intractably large state and action spaces. Although the artificial intelligence and operations research communities have independently proposed two successful frameworks for solving dynamic stochastic optimization problems---Monte Carlo tree search (MCTS) and mathematical optimization (MO), respectively---the relative merits of these two approaches are not well understood. In this paper, we adapt both MCTS and MO to a problem inspired by tactical wildfire and management and undertake an extensive computational study comparing the two methods on large scale instances in terms of both the state and the action spaces. We show that both methods are able to greatly improve on a baseline, problem-specific heuristic. On smaller instances, the MCTS and MO approaches perform comparably, but the MO approach outperforms MCTS as the size of the problem increases for a fixed computational budget.  \end{abstract}


\end{frontmatter}

\section{Introduction}
\label{sec:Intro}
Dynamic resource allocation (DRA) problems are problems where one must assign resources to tasks or requests over some finite time horizon. Many important real-world problems can be cast as DRA problems, including applications in air traffic control \cite{barnhart2003applications, Bertsimas1998}, network revenue management (see, e.g., \cite{talluri2004theory} for an overview), scheduling \cite{Bertsimas2013, doi:10.1080/0740817X.2010.504690}  and logistics, transportation  and fulfillment \cite{acimovic2012making, powell2002adaptive}. DRA problems are notoriously difficult to solve for two reasons. First, many real-world DRA problems exhibit stochasticity, i.e., the requests to be processed may arrive randomly according to some stochastic process that, itself, depends on where resources are allocated.  Second, many real-world DRA problems exhibit extremely large state and action spaces, making solution by traditional dynamic programming methods infeasible \cite{Bellman1957,Powell2011}.  A number of scientific communities, particularly within artificial intelligence and operations research, have sought more sophisticated techniques for addressing DRA and other dynamic stochastic optimization problems.  

Within the AI community, one approach for dynamic stochastic optimization problems that has received increasing attention in the last 15 years is a method known as Monte Carlo tree search (MCTS) \cite{Coulom2007, Browne2012}. In any dynamic stochastic optimization problem, one can represent the possible trajectories of the system---the state at each decision epoch and the actions taken at those epochs---as a tree, where the root represents the initial state. In MCTS, one iteratively builds an approximation to this tree and uses it to inform the choice of action.
In our opinion, MCTS's effectiveness stems from two key features: 1) ``bandit upper confidence bounds" (see \cite{Auer2002, Kocsis2006}) can be used to balance exploration and exploitation in learning, 2) application-specific heuristics and knowledge can be used to customize the base algorithm \cite{Browne2012}.  
Moreover, the MCTS algorithm is very flexible and can easily be tailored to a variety of problems.  Indeed, the only technical requirement for implementing MCTS is a generative model that, given a state and an action at a given decision epoch, generates a new state for the next epoch and an immediate reward received.  This flexibility makes MCTS particularly attractive as a general purpose methodology. 

Most importantly, MCTS has been extremely successful in a number of applications, particularly in designing expert computer players for difficult games such as Go \cite{Gelly2011,Lee2009,Enzenberger2010}, Hex \cite{Arneson2009,Arneson2010}, Kriegspiel \cite{Ciancarini2010}, and Poker \cite{Rubin2011}.  Although, we think it is fair to say that MCTS is the one of the top performing, general purpose algorithms for this class of games, we observe that games like Go and Hex are qualitatively very different from DRAs.  Namely, unlike typical DRA problems, the state of these games does not evolve stochastically, and, furthermore, the size of the feasible action space is often much smaller than for typical DRA problems. For example, in the instances of Go studied in \cite{Gelly2011}, the action branching factor is at most 81, whereas in the DRA instance we consider, a typical branching factor is approximately 230 million (cf. Eq.~\eqref{eq:stirling}).  This raises a question about the performance of MCTS for this new class of real-world problems.  

On the other hand, within the operations research community, the study of DRAs has proceeded along different lines.  A prominent stream of research is based upon mathematical optimization (MO).  In contrast to MCTS which only requires access to a generative model of the stochastic system, MO approaches model the dynamics of the system \emph{explicitly} via a constrained optimization problem.  The solution to this optimization problem then yields a control policy for the system.  This paradigm has a long history within the dynamic control literature (see, e.g., \cite{Bertsekas1995} for an overview).  

In the special case of stochastic, dynamic resource allocation problems, there a number of competing proposals to incorporate the stochastic evolution including robust optimization / control \cite{BertsimasBrownCaramanis2011, ben2009robust, nilim2005robust, grieder2003robust} and chance constrained optimization \cite{charnes1959chance}.  In what follows, however, we focus on a third proposal, sometimes called \emph{model predictive control} \cite{bemporad2006model}.  Specifically, we replace uncertain parameters in a MO formulation with their expected values and periodically re-solve the optimization problem for an updated policy as the true system evolves.
Variants of this paradigm are known by many names in the literature including fluid approximation \cite{Gallego1994,Avram1995}, certainty equivalent control (\cite[Chapt. 6]{Bertsekas1995}).  
These (re-)optimization frameworks are well-known to have excellent practical performance in application domains like queueing \cite{chen1993dynamic, shah2012switched} and network revenue management \cite{Ciocan2012}, and in some special cases, additionally enjoy strong theoretical guarantees on performance (e.g., \cite{Ciocan2012, Gallego1994, jasin2012re, maglaras2006dynamic, reiman2008asymptotically, Talluri1998}).  

In any case, regardless of the particular proposal for incorporating stochasticity, the widespread use and success of MO approaches for DRAs contrasts strongly with a lack of computational experience with MCTS for DRAs.

In this paper, we aim to understand the relative merits of both the above MCTS and MO approaches by applying them to a challenging DRA problem inspired by tactical wildfire management.  
Our interest in tactical wildfire management as a benchmark problem stems simultaneously from practical importance of this application area and the computational difficulty of the underlying DRA.  
Indeed, the severity of wildland fires has been steadily increasing in recent years, particularly in the southwestern US \cite{gorte2013}, and as a result, US federal government spending on wildfire management has also been increasing, amounting to \$3.5 billion in 2013 \cite{Bracmort2013}. Suppression costs are only part of the total cost of wildfire; combined costs of loss of property, damage to the environment and loss of human life are estimated to be many times larger \cite{TrueCostWildfire2010}.  Consequently, there has been renewed interest in models for realtime decision support tools to assist in tactical wildfire management, namely how best to utilize fire suppression resources to contain and extinguish a fire.  Unfortunately, this is an extremely challenging problem both from the complexity of the system dynamics and uncertainty in the underlying data.  

Although there have been a number of empirically validated, deterministic models for wildfire spread proposed (e.g., \cite{Finney2004, Prometheus2009}), there have been fewer works that incorporate the stochastic elements of fire spread \citep{Boychuck2008,doi:10.1139/x2012-032,Fried2006}. Most works focus on developing models for the rate and simulation of the spread of wildfire.  Fewer consider the associated problem of managing suppression resources.  A notable exception is the research stream \cite{Hu:2009:ISO:1596519.1596524, 8582580020130201}, which proposes 1) an optimization formulation of the \emph{initial attack} problem, i.e., the problem of determining how many and what kind of suppression resources to allocate to an ongoing fire and 2) a simulation engine to study how specific, user-specified heuristic suppression rules affect the growth of a fire.  To the best of our knowledge, the authors do not address the tactical problem of \emph{optimally} controlling suppression resources as the fire evolves.  In other words, the authors do not consider the underlying DRA.  

By contrast, we develop a model for the underlying DRA as a particular Markov decision problem (MDP) which incorporates several of the most prominent features of existing wildfire spread models.  We then compare the performance of MCTS and MO inspired algorithms for this problem.  We stress that this is a particularly challenging MDP, having many possible states and many possible actions for even small instances of the problem.  Our computational experience with this benchmark both elucidates some of the features of these algorithms and illustrates how they might be incorporated into a larger decision support system like \cite{Hu01032012}.  

%

We summarize our contributions as follows:
\begin{enumerate}
	\item We propose a flexible MDP to model tactical wildfire management that captures the most important deterministic and stochastic features of wildfire spread.  We also propose a simple, high-quality, and customized heuristic for this problem that may be of independent interest to practitioners.  
	\item We develop an MCTS-based approach for the above MDP.  To the best of our knowledge, this represents the first application of MCTS to a DRA problem motivated by a real-world application.  Towards this end, we combine a number of classical features of MCTS, such as bandit upper confidence bounds, with new features such as double progressive widening \cite{Couetoux2011}.  We also propose a novel action generation approach to cope with the size of the state and action spaces of the DRA.  
	\item We propose a mathematical optimization formulation that approximates the original discrete and stochastic elements of the MDP by suitable continuous and deterministic counterparts.  Although this approximation is in the same spirit as other fluid approximation literature in operations research (\cite{Avram1995, Gallego1994}), our particular formulation incorporates elements of a linear dynamical system which we believe may be of independent interest in other DRA problems.  
	\item Through extensive computational experiments, we show the following:
		\begin{enumerate}
			\item The MCTS and MO approaches both produce high-quality solutions, generally performing as well or better than our customized heuristic for this problem. MCTS and MO perform comparably when the problem instance is small.  With a fixed computational budget, however, the MO approach begins to outperform the MCTS approach as the size of the problem instance grows, either in state space or action space. Indeed, for very large action spaces, MCTS can begin to perform worse than our baseline heuristic.  The MO approach, by comparison, still performs quite well. 
			\item The choice of hyperparameters in the MCTS algorithm---such as the exploration bonus and the progressive widening parameters---can significantly affect overall performance of the algorithm.  The interdependence between these parameters is complex, and they cannot, in general, be selected independently.  Some care must be taken to appropriately tuning the algorithm to a specific DRA problem.  
		\end{enumerate}
\end{enumerate}

The rest of this paper is organized as follows. In Section~\ref{sec:Dynamics}, we introduce our MDP formulation of the tactical wildfire management problem that we use to compare the MCTS and MO approaches. In Sections~\ref{sec:MCTS} and \ref{sec:Optimization}, we describe the MCTS and MO approaches, respectively. In Section~\ref{sec:Numerics}, we describe our computational experiments and report on the results of these experiments. Finally, in Section~\ref{sec:Conclusion}, we summarize the main contributions of the work and highlight promising directions for future research.


\section{Tactical Wildfire Dynamics as Dynamic Resource Allocation}
\label{sec:Dynamics}
\usetikzlibrary{arrows,shapes,automata}
\tikzstyle{state}=[minimum width=12mm,circle,draw=black]
\tikzstyle{action}=[minimum width=12mm,minimum height=12mm,rectangle,draw=black]
\tikzstyle{reward}=[minimum width=16mm,minimum height=16mm,diamond,draw=black]

\label{sec:dynamicmodel}

In what follows, we consider a simplified model in order to more clearly highlight the difference between the MCTS and MO approaches in our experiments while still capturing the key features of wildfire propagation, e.g., stochastic evolution, wind and topography effects and dependence on fuel.  Although we do not explore the ideas here, it would be (at least conceptually) straightforward to extend our approach to a higher fidelity model incorporating moisture content and surface fuel-type as in \cite{Rothermel1972} by using techniques similar to those developed in \cite{Hu01032012}.  

Specifically inspired by the stochastic model of \citep{Boychuck2008} for wildland fire simulation, we model tactical wildland fire management as a Markov decision process. 
We partition the landscape into a grid of cells $\mathcal{X}$. There are two attributes for each cell $x \in \mathcal{X}$:
\begin{itemize}
\item $B(x)$, a Boolean variable indicating whether the cell is currently burning, and
\item $F(x)$, an integer variable indicating how much fuel is remaining in the cell.
\end{itemize}
The collection of these two attributes over all cells in the grid represents the state of the Markov decision process.  

We further assume a set of suppression resources $\mathcal{I}$ fighting the wildland fire.  To simplify the model, we will treat all suppression resources as identical, but it is straightforward to extend to the case heterogenous resources with different capabilities as in \citep{MARTIN-FERNANDEZ2002, Hu01032012}.
The decisions of our MDP correspond to assigning resources to cells.  For each $i \in \mathcal{I}$, let $a^{(i)} \in \mathcal{X}$ denote the cell to which we assign suppression resource $i$.  We assume that any resource can be assigned to any cell at any time step, i.e., that the travel time between cells is negligibly small compared to the decision interval of the MDP.

Once ignited, a cell consumes fuel at a constant rate.  Once the fuel is exhausted, the cell extinguishes.  Since fuel consumption occurs a a constant rate, without loss of generality, we can rescale the units of time to make this rate equal to unity.  Thus, we model the evolution of fuel in the model by
\begin{equation}
F_{t+1}(x) = \begin{cases}
	F_t(x) & \mbox{if } \neg B_t(x) \vee F_t(x) = 0\\
	F_t(x) - 1 & \mbox{otherwise.}
\end{cases}
\end{equation}
Notice this evolution is deterministic given $B_t(x)$.  

\begin{figure}[!ht]
\begin{center}
\begin{tikzpicture}[font=\footnotesize, node distance=20mm,->]

\node [state] (T1) {false};
\node [state,right of=T1] (T2){true};

\draw (T1) edge [bend left] node [above] {$\rho_1$} (T2);
\draw (T1) edge [loop above] node [above] {$1-\rho_1$} (T1);

\draw (T2) edge [bend left] node [above] {$\rho_2$} (T1);
\draw (T2) edge [loop above] node [above] {$1-\rho_2$} (T2);

\end{tikzpicture}
\end{center}
\caption{$B(x)$ transition model.\label{fig:Btran}}
\end{figure}
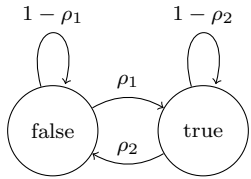

The evolution of $B_t(x)$, however, is stochastic.  Figure~\ref{fig:Btran} shows the probabilistic transition model for $B_t(x)$ where
\begin{equation}
\rho_1 = \begin{cases}
	1 - \prod_{y} (1-P(x,y) B_{t}(y)) & \mbox{if } F_t(x) > 0 \\
	0 & \mbox{otherwise}
\end{cases}
\end{equation}
and
\begin{equation}
\rho_2 = \begin{cases}
	1 & \mbox{if } F_t(x) = 0 \\
	1 - \prod_{i} (1-Q(x)\delta_x ( a^{(i)} ) ) & \mbox{otherwise.}
\end{cases}
\end{equation}
Here, $P(x,y)$ is the probability that a fire in cell $y$ ignites a fire in cell $x$. Generally, only the neighbors of $x$ can ignite $x$, and so we expect $P(x,y)$ to be sparse.  The specification of $P(x,y)$ can capture the tendency of a fire to propagate primarily in one direction due to wind or sloping terrain.  $Q(x)$ is the probability that a suppression effort on cell $x$ successfully extinguishes the cell.  We assume that the probability of success for multiple attempts on the same cell are independent.  

We stress that under these dynamics, cells that have been previously extinguished by a suppression team may later reignite.  

The reward for a cell burning is $R(x)$ (always negative) and the total reward received every step is $\sum_{t, x} B_t(x)R(x)$.  We can vary the reward across the grid to represent a higher cost of a fire in particular areas. For example, we may penalize a fire in a populated area more heavily.


\section{Monte Carlo Tree Search}
\label{sec:MCTS}
Our review of MCTS is necessarily brief.  Please see \cite{Browne2012} for a more thorough survey.  In particular, we will focus on a variation of MCTS that employs \emph{double progressive widening}, a technique that explicitly controls the branching factor of the search tree ~\cite{Couetoux2011}. This variation is specifically necessary when the action space is continuous or so large that all actions cannot possibly be explored.  This is frequently the case in most DRA problems.  We note that there exist other strategies for managing the branching factor, such as the pruning strategy employed in the ``bandit algorithm for smooth trees'' (BAST) method of \cite{Coquelin2007}; we do not consider these here. 

We first describe MCTS with double progressive widening in general, and then discuss two particular modifications we have made to the algorithm to tailor it to DRAs.  

\subsection{MCTS with double progressive widening}
Algorithm~\ref{alg:mcts-dpw} involves running many simulations from the current state while updating an estimate of the state-action value function $Q(s, a)$. We use a generative model $G$ to produce samples of the next state $s'$ and reward $r$ given the current state $s$ and action $a$. We draw samples $(s', r) \sim G(s,a)$. All of the information about the state transitions and rewards is represented by $G$; the state transition probabilities and expected reward function are not used directly. There are three stages in each simulation: search, expansion, and rollout.
\begin{algorithm}[!ht]
\caption{\label{alg:mcts-dpw}Monte Carlo tree search with double progressive widening}
\begin{algorithmic}
\Function{MonteCarloTreeSearch}{$s, d$}
\Loop
	\State $\Call{Simulate}{s, d}$
\EndLoop
\State \Return $\argmax_a Q(s, a)$
\EndFunction

\Function{Simulate}{$s, d$}
\If{$d = 0$}
	\State \Return $0$
\EndIf	
\If{$s \not\in T$}
	\State $T = T \cup \{s\}$
	\State $N(s) \gets N_0(s)$
	\State \Return $\Call{Rollout}{s, d}$	
\EndIf

\State $N(s) \gets N(s) + 1$
\If{$\| A(s) \| < k N(s)^\alpha$}
	\State $a \gets \Call{Getnext}{s, Q}$
	\State $(N(s, a), Q(s, a), V(s, a)) \gets (N_0(s, a), Q_0(s, a), V(s, a))$
	\State $A(s) = A(s) \cup \{a\}$
\EndIf

\State $a \gets \argmax_a Q(s, a) + c \sqrt{\frac{\log N(s)}{N(s,a)}}$
\If{$ \| V(s,a) \| < k' N(s,a)^{\alpha'}$}
	\State $(s', r) \sim G(s, a)$ \label{line:generative mcts}
	\If{$s' \not\in V(s,a)$}	
		\State $V(s,a) = V(s,a) \cup \{s'\}$
		\State $R(s,a,s') \gets r$
		\State $N(s,a,s') \gets N_0(s,a,s')$
	\Else
		\State $N(s,a,s') \gets N(s,a,s')+1$
	\EndIf
\Else
	\State $s' \gets \Call{Sample}{N(s,a,\cdot)}$
	\State $r \gets R(s,a,s')$
	\State $N(s,a,s')  \gets N(s,a,s')+1$
\EndIf

\State $q \gets r + \gamma \Call{Simulate}{s', d - 1}$
\State $ N(s,a) \gets  N(s, a) + 1$
\State $Q(s, a) \gets Q(s, a) + \frac{q - Q(s, a)}{N(s, a)}$
\State \Return $q$
\EndFunction

\Function{Rollout}{$s, d$}
\If{$d = 0$}
	\State \Return $0$
\EndIf
\State $a \sim \pi_0(s)$
\State $(s', r) \sim G(s, a)$ \label{line:generative mcts2}
\State \Return $r + \gamma \Call{Rollout}{s', a, d - 1}$
\EndFunction
\end{algorithmic}
\end{algorithm}

\subsubsection{Search}
If the current state in the simulation is in the set $T$ (initially empty), we enter the search stage. Otherwise, we proceed to the expansion stage. During the search stage, we update $Q(s, a)$ for the states and actions visited and tried in our search. We also keep track of the number of times we have visited a state $N(s)$ and the number of times we have taken an action from a state $N(s, a)$. 

During the search, the first progressive widening controls the number of actions considered from a state. To do this, we generate a new action if $\| A(s) \| < k N(s)^\alpha$, where $k$ and $\alpha$ are parameters that control the number of actions considered from the current state and $A(s)$ is the set of actions tried from $s$. When generating a new action, we add  it to the set $A(s)$, and initialize $N(s, a)$ and $Q(s, a)$ with $N_0(s, a)$ and $Q_0(s, a)$, respectively. The functions $N_0$ and $Q_0$ can be based on prior expert knowledge of the problem; if none is available, then they can both be initialized to $0$. We also initialize the empty set $V(s,a)$, which contains the set of states $s'$ transitioned to from $s$ when selecting action $a$. A default strategy for generating new actions is to randomly sample from candidate actions. If $\| A(s) \| \ge k N(s)^\alpha$, then we execute the action that maximizes
\begin{equation*}
Q(s, a) + c \sqrt{\frac{\log N(s)}{N(s,a)}}\,, 
\end{equation*}
where $c$ is a parameter that controls the amount of exploration in the search. The second term is an \emph{exploration bonus} that encourages selecting actions that have not been tried as frequently.

Next, we draw a sample $(s', r) \sim G(s,a)$, if $ \| V(s,a) \| < k' N(s,a)^{\alpha'}$. In this second progressive widening step, the parameters $k'$ and $\alpha'$ control the number of states transitioned to from $s$. If $s'$ is not a member of $V(s,a)$, we add it to the set $V(s,a)$, initialize $R(s,a,s')$ to $r$, and initialize $N(s,a,s')$ with $N_0(s,a,s')$.  If $s'$ is a member of $V(s,a)$, then we increment $N(s,a,s')$. However, if $\| V(s,a) \| \ge k' N(s,a)^{\alpha'}$, then we select $s'$ from $V(s,a)$ proportional to $N(s,a,s')$.

\subsubsection{Expansion}
Once we have reached a state that is not in the set $T$, we generate a new action $a$ available from the state, add it to the set $A(s)$, and initialize $N(s, a)$ and $Q(s, a)$ with $N_0(s, a)$ and $Q_0(s, a)$, respectively. We then add the current state to the set $T$.

\subsubsection{Rollout}
After the expansion stage, we simply select actions according to some \emph{rollout} (or default) policy $\pi_0$ until the desired depth is reached. Typically, rollout policies are stochastic, and so the action to execute is sampled $a \sim \pi_0(s)$. The rollout policy does not have to be close to optimal, but it is a way for an expert to bias the search into areas that are promising. The expected value is returned and is used in the search to update the value for $Q(s, a)$ used by the search phase.

Simulations are run until some stopping criteria is met, often simply a fixed number of iterations. We then execute the action that maximizes $Q(s, a)$. Once that action has been executed, we can rerun Monte Carlo tree search to select the next action. It is common to carry over the values of $N(s, a)$, $N(s)$, and $Q(s, a)$ computed in the previous step.

\subsection{Tailoring MCTS to DRAs}
We next discuss two modifications we have found to be critical when applying MCTS to DRAs.  

\subsection{Action Generation}
As previously mentioned, the default strategy for generating new actions during the search stage of MCTS involves randomly sampling an action from all candidate actions.  In DRAs where the action space may be very large, this strategy is inefficient; we may need to search many actions before identifying a high-quality choice.  Rather, we would like to bias the sampling towards potential actions that we believe may perform well.  One way to identify such actions is via some application specific heuristic.  

We follow a slightly different approach.  Specifically, consider MCTS after several iterations.  The current values of $Q(s, a)$ 
provide a (noisy) estimate of the value function, and, hence, can be used to approximately identify promising actions.  
Consequently, we use these estimates to bias our sampling procedure through a sampling scheme inspired by genetic algorithm search heuristics~\citep{raey}. Our strategy, described in Algorithm~\ref{alg:getnextga}, involves generating actions using one of three approaches: with probability $u'$ an existing action in the search tree is mutated, with probability $u''$ two existing actions in the search tree are recombined, or a new action is generated from the default strategy. Mutating involves randomly changing the allocation of one or resources. Recombining involves selecting a subset of allocations from two actions and combining the two subsets. When mutating or recombining, we select the existing action (or actions) from $A(s)$ using tournament select where the fitness for each action is proportional to $Q(s,a)$. Note that it is permissible to use other methods such as softmax, too.
\begin{algorithm}[!ht]
\caption{\label{alg:getnextga}Action Generation}
\begin{algorithmic}
\Function{Getnext}{$s, Q$}
\State $u \sim U(0,1)$
\If{$u < u'$}
	\State $a' \gets \Call{Sample}{Q(s,\cdot)}$
	\State $a \gets \Call{Mutate}{a'}$
\ElsIf {$u < u'+u''$}
	\State $a' \gets \Call{Sample}{Q(s,\cdot)}$
	\State $a'' \gets \Call{Sample}{Q(s,\cdot)}$
	\State $a \gets \Call{Recombine}{a', a''}$
\Else
	\State $a \sim \pi_0(s)$
\EndIf
\State \Return $a$
\EndFunction
\end{algorithmic}
\end{algorithm}

Our numerical experiments confirm that our proposed action generation approach significantly outperform the default strategies in our benchmark problem.

\subsection{Rollout Policy}
\label{subsec:dp_rollout}
In many papers treating MCTS, it is argued that even if the heuristic used for the rollout policy is highly suboptimal, given enough time the algorithm will converge to the correct state-action value function $Q(s, a)$.  In DRAs with combinatorial structure (and, hence, huge state  and action spaces), it may take an extremely long time for this convergence.  (Similar behavior was observed in \cite{coquelin2007bandit}).  Indeed, we have observed that for DRAs that having a good initial rollout policy makes a material difference in the performance of MCTS.  Unfortunately, designing a good heuristic seems to be an application specific task.  

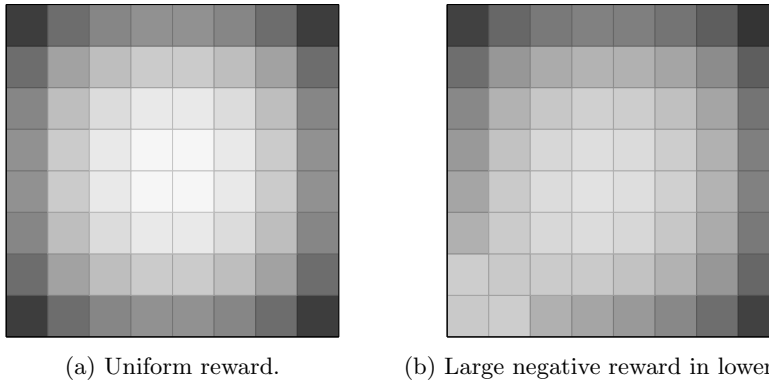
\begin{figure}[!ht]
\centering
\begin{subfigure}{0.46\textwidth}
\centering
	\begin{tikzpicture}
	\begin{axis}[width=6cm, height=6cm,view={0}{90},xtick=\empty,ytick=\empty,colormap/blackwhite]
		\addplot3[surf] file {fw-heuristic-1.dat};
	\end{axis}
	\end{tikzpicture}
	\caption{Uniform reward.}
\end{subfigure}
\hspace{0.05cm}
\begin{subfigure}{0.46\textwidth}
\centering
	\begin{tikzpicture}
	\begin{axis}[width=6cm, height=6cm,view={0}{90},xtick=\empty,ytick=\empty,colormap/blackwhite]
		\addplot3[surf] file {fw-heuristic-3.dat};
	\end{axis}
	\end{tikzpicture}
	\caption{Large negative reward in lower left.}
\end{subfigure}
\caption[Example heuristic weights.]{Example heuristic weights. Lighter cells correspond to higher weights.\label{fig:fw-heuristic}}
\end{figure}

For our benchmark, we use a heuristic that involves assigning a weight to each cell $x$ 
\begin{equation}
W(x) = \sum_y \frac{R(y)}{D(x,y)} \eqpunct{,}
\end{equation}
where $D(x,y)$ is the shortest path between $x$ and $y$ assuming that the distance between adjacent cells is $P(x,y)$. We compute the values $D(x,y)$ offline using a graph analysis algorithm, such as the Floyd-Warshall algorithm~\cite{Floyd1962}, and consequently term this heuristic the FW heuristic. Figure~\ref{fig:fw-heuristic} shows example weights assigned to different cells on an eight by eight grid with two different reward profiles. The heuristic performs generally well because it prioritizes allocating resources to cells that are near large negative reward cells (i.e., populated areas).

The rollout policy involves selecting the cells that are burning and assigning resources to the highest weighted cells. We are also able to randomly sample from the weights to generate candidate actions. 

In our experiments, we have observed that the heuristic performs fairly well, and, consequently, may be of independent interest to the wildfire suppression community.  


\section{A Mathematical Optimization Approach}
\label{sec:Optimization}
In this section, we present an optimization-based solution approach for the tactical wildland fire management problem. In this approach, we formulate a deterministic optimization problem approximating the original MDP.  At each decision step, we re-solve this approximation based on the current state of the process and select the first prescribed allocation.  

The key feature of the formulation is the use of a deterministic, ``smoothed'' version of the dynamics presented in Section~\ref{sec:Dynamics}.  Rather than modeling the state with a discrete level of fuel and binary state of burning, we model fuel as a continuous quantity and model a new (continuous) \emph{intensity} level of each cell representing the rate at which fuel is consumed. Other authors have used similar ideas when motivating various fluid approximations in the operations research literature. For example, continuous fluid approximations have been used to study the control of queuing networks \cite{Avram1995}, where the size of each queue in the network can be modeled continuously through systems of differential equations, and the decision at each time for each server in the network is the ``rate'' at which each queue is being served/emptied. Another example comes from revenue management (RM). In a typical RM problem, the decision maker needs to sell some fixed inventory of a product, and must at each point in time set a price for this product. The price determines the rate of a stochastic demand process, with the goal of maximizing the total revenue at the end of the selling period. The optimal solution of the exact problem is, with few exceptions, extremely difficult to obtain. However, if the dynamics of the problem are relaxed so that the demand is continuous and arrives at a deterministic rate that varies with the price, one can obtain useful bounds and near-optimal pricing policies for the exact, \emph{stochastic} problem \cite{Gallego1994}.

Most importantly, smoothing the dynamics allows for two important simplifications.  First, we can replace the probabilistic dynamics in Section~\ref{sec:Dynamics} with simpler, deterministic dynamics governing the evolution of the intensity.  Second, and most importantly, we no longer have to consider the entire exponentially large state space of the MDP, but rather only its evolution along the one path determined by the deterministic dynamics.

%

\subsection{Optimization Model}
\label{subsec:optimization_deterministic}
Let $A_t(x,i)$ be a binary variable that is 1 if suppression resource $i \in \Ical$ is assigned to cell $x$ at time $t$ and 0 otherwise; this is the main decision variable of the problem. Recall that $F_t(x)$ denotes the amount of fuel available at the start of period $t$ in cell $x$.  Furthermore, let 
$I_t(x)$ represent the intensity of the fire in cell $x \in \Xcal$ at time $t$. Intensity is a continuous decision variable that will be determined by the optimization algorithm.  Unlike the original MDP formulation, it is no longer possible to rescale the parameters so that exactly one unit of fuel is consumed per period.  (We discuss how to calibrate the fuel appropriately next.)  

Some additional notation is required to describe the evolution of $I_t(x)$.  
Define $\Ncal(x)$ as the set of cells that are neighbors of $x$ in the sense of fire transmission, i.e. $\Ncal(x) = \{ y : P(x, y) > 0 \}$.  Let $\zeta_t(y,x) \in [0,1]$ be the rate at which the intensity $I_{t-1}(y)$ of cell $y$ at time $t-1$ contributes to the intensity $I_t(x)$ of cell $x$ at time $t$.
Furthermore, let $\tilde{\zeta}_t(x,i) \in [0,1]$ be the relative reduction in intensity when suppression team $i$ is assigned to cell $x$ at time $t$.  
Finally, let $I_0(x) = 1$ if cell $x$ is burning at time $0$, and $I_0(x) = 0$ if cell $x$ is not burning.


With this notation, our formulation is
\begin{subequations}
\begin{align}
& \text{minimize} & & \sum_{x \in \Xcal} \sum_{t=0}^T R(x) I_t(x) \label{prob:intensity_obj}\\
& \text{subject to} & & I_t(x) \geq  I_{t-1}(x) + \sum_{y \in \Ncal(x)} I_{t-1}(y) \cdot \zeta_t(y,x) \nonumber \\
& & &  - \sum_{i \in \Ical} A_{t-1}(x,i) \cdot \bar{I}_t(x) \cdot \tilde{\zeta}_t(x,i) - \left(F_0(x) + \sum_{y \in \Ncal(x)} F_0(y)\right) \cdot z_{t-1}(x) , \nonumber \\ 
& & &  \quad \forall x \in \Xcal, \ t \in \{ 1, \dots, T\}, \label{prob:intensity_intensitydynamics} \\[0.25em]
& & & F_t(x) = F_0(x) - \sum_{t'=0}^{t-1} I_{t'}(x), \quad \forall x \in \Xcal, \ t \in \{0, \dots, T\}, \label{prob:intensity_fueldefinition} \\
& & & F_t(x) \geq \delta \cdot ( 1 - z_t(x)), \quad \forall x \in \Xcal, \ t \in \{0, \dots, T\}, \label{prob:intensity_fuelzclamp1} \\
& & & F_t(x) \leq \delta \cdot z_t(x) + F_0(x) \cdot ( 1 - z_t(x) ), \quad \forall x \in \Xcal, \ t \in \{0, \dots, T\}, \label{prob:intensity_fuelzclamp2} \\
& & & I_{t+1}(x) \leq F_0(x) \cdot (1 - z_t(x)), \quad \forall x \in \Xcal, \ t \in \{0, \dots, T\}, \label{prob:intensity_fuelintensityforcing}\\
& & & \sum_{x \in \Xcal} A_t(x,i) \leq 1, \quad \forall t \in \{0, \dots, T\}, \ i \in \Ical, \label{prob:intensity_assignmentofteams}\\
& & & I_t(x),\ F_t(x) \geq 0, \quad \forall x \in \Xcal, \ t \in \{0, \dots, T\}, \\
& & & z_t(x) \in \{0,1\}, \quad \forall x \in \Xcal,\ t \in \{0, \dots, T\}, \\
& & & A_t(x,i) \in \{0,1\}, \quad \forall x \in \Xcal,\ i \in \Ical,\ t \in \{0, \dots, T\}.
\end{align}
\label{prob:intensity}
\end{subequations}

Here, $\delta > 0$ is a small threshold chosen so that a cell with less than $\delta$ units of fuel cannot burn. Consequently, the binary decision variable $z_t(x)$ represents the indicator that the fuel at time $t$ in cell $x$ is below $\delta$.  Finally in the spirit of ``big-M'' constraints, $\bar{I}_t(x)$ is an upperbound on the maximal value attainable by $I_t(x)$.  We will discuss how to compute this value shortly.

The constraints have the following meaning:
\begin{itemize}
\item Constraint~\eqref{prob:intensity_intensitydynamics} expresses the one-step dynamics of the fire intensity in region $x$ at time $t$.  
Although we have written the constraint in inequality form, it is not difficult to see that in an optimal solution, the constraint will always be satisfied at equality since the objective is a sum of the intensities over all periods and regions weighted by the (positive) importance factors. 

The first two terms of the right-hand side represent that---without intervention and without regard for fuel---the intensity of a cell one step into the future is the current intensity ($I_{t-1}(x)$) plus the sum of the intensities of the neighboring cells weighted by the transmission rates ($\sum_{y \in \Ncal(x)} I_{t-1}(y) \cdot \zeta_t(y,x)$). If suppression team $i$ is assigned to cell $x$ at $t-1$, then $A_t(x,i) = 1$ and the intensity is reduced by $\tilde{\zeta}_t(x,i) \cdot \bar{I}_t(x)$. If the cell's fuel is below $\delta$ at time $t-1$, then $z_t(x) = 1$, and the intensity is reduced by $F_0(x) + \sum_{y \in \Ncal(x)} F_0(y)$; since the intensity of a cell is upper bounded by the initial fuel of that cell, the term $- \left(F_0(x) + \sum_{y \in \Ncal(x)} F_0(y)\right) \cdot z_{t-1}(x)$ ensures that whenever the fuel $F_t(x)$ drops below $\delta$, this constraint becomes vacuous. 
\item Constraint~\eqref{prob:intensity_fueldefinition} is the equation for the remaining fuel at a particular time point as a function of the intensities (intensity is assumed to be the fuel burned in a particular time period).
\item Constraint~\eqref{prob:intensity_fuelzclamp1} and \eqref{prob:intensity_fuelzclamp2} are forcing constraints that force $F_t(x)$ to be between $\delta$ and $F_0(x)$ if $z(t) = 0$, and between $0$ and $\delta$ if $z(t) = 1$.
\item Constraint~\eqref{prob:intensity_fuelintensityforcing} ensures that if there is insufficient fuel in cell $x$ at period $t$, then the intensity at that cell in the next time point is zero. If there is sufficient fuel, then the constraint is vacuous (the intensity is at most $F_0(x)$, which is already implied in the formulation).
\item Constraint~\eqref{prob:intensity_assignmentofteams} ensures that each team in each period is assigned to at most one cell.
\item The remaining constraints ensure that the fuel and intensity are continuous nonnegative variables, and that the sufficient fuel and team assignment variables are binary. 
\end{itemize}
The objective \eqref{prob:intensity_obj} is the sum of the intensities over all of the time periods and over all cells, weighted by the importance factor of each cell in each time period. 

Problem~\eqref{prob:intensity} is a mixed integer linear optimization model with two sets of binary variables: the $A_t(x,i)$ variables, which model the assignment of suppression teams to cells over time, and the $z_t(x)$ variables, which model the loss of fuel over time. 
Although mixed linear optimization is 
not solvable in polynomial time, there exist high-quality open-source and commercial solvers which are able to solve such problems extremely efficiently in practice, even for very large instances.  

In highly resource constrained environments when it is not possible to solve the above model to optimality, we can still obtain extremely good approximate solutions by relaxing the $A_t(x,i)$ variables to be continuous within the unit interval $[0,1]$. Then, given an optimal solution with fractional values for the $A_t(x,i)$ variables at $t = 0$, we can compute a score $v(x)$ for each cell $x$ as $v(x) = \sum_{i \in \Ical} A_0(x,i)$. We then assign suppression teams to the $|\Ical|$ cells with the highest values of the index $v$.  Indeed, we will follow this strategy in Section~\ref{sec:Numerics}.  

\subsection{Calibrating the Model}
Given parameters for the original MDP formulation of the tactical wildland fire management problem, parameters for our nominal optimization formulation are obtained as follows:
\begin{itemize}
\item $\bar{I}_t(x)$ is computed by iterating a modified version of the one-step recursion, assuming that there is no intervention and infinite fuel:
\begin{equation*}
\bar{I}_t(x) = \bar{I}_{t-1}(x) + \sum_{y \in \Ncal(x)} \bar{I}_{t-1}(y),
\end{equation*}
where $\bar{I}_0(x) = I_0(x)$. In this modified one-step recursion, each transmission rate $\zeta_t(y,x)$ is essentially assumed to be 1, which is the highest value it can be. 
\item $F_0(x)$ is obtained by summing the fuel threshold $\delta$ and the $\bar{I}_t(x)$ values over the horizon $t = 0,1, \dots, \min\{T, F(x)\}$, where $F(x)$ is the number of periods that cell $x$ can burn into the future according to the original MDP dynamics:
\begin{equation*}
F_0(x) = \delta + \sum_{t=0}^{\min\{T,F(x)\} } \bar{I}_t(x).
\end{equation*}
Intuitively, since the intensity $I_t(x)$ can be thought of as how much fuel was consumed by the fire in cell $x$ at time $t$, the initial fuel value $F_0(x)$ can be thought of as a limit on the cumulative intensity in a cell over the entire horizon. Once the cumulative intensity has reached $\sum_{t=0}^{\min\{T,F(x)\} } \bar{I}_t(x)$, the fuel in the cell enters the interval $[0,\delta]$, at which point the variable $z_t(x)$ is forced to 1 and the intensity is forced to zero for all remaining time periods. 
\item $\zeta_t(y,x)$ is set to $P(x,y)$ (the transmission probability from $y$ to $x$) for each $t$.
\item $\tilde{\zeta}_t(x,i)$ is set to $Q(x)$ (the probability of successful extinguishing a fire in cell $x$) for each period $t$ and each suppression team $i$.
\item $\delta$ is set to $0.1$.
\end{itemize}



\section{Numerical Comparisons}
\label{sec:Numerics}
This section presents experiments comparing Monte Carlo tree search (MCTS) and the mathematical optimization formulation (MO).  We seek to understand their relative strengths and weaknesses as well as how the user-defined parameters of each approach, such as the exploration bonus $c$ for MCTS, affect the performance of the algorithm.  Before proceeding to the details of our experiments, we summarize our main insights here:
\begin{itemize}
	\item Overall, MO performs as well as or better than MCTS. For even moderate computational budgets, MO generates high quality solutions. 
	\item Although the MCTS approach works well for certain smaller examples, its performance can degrade for larger examples (with a fixed, computational budget).  Moreover, the dependence on the algorithm on its underlying hyperparameters is complex.  The optimal choice of the exploration bonus and progressive widening factors depends both on the underlying heuristic used in the rollout as well as available computational budget.
\end{itemize}

\subsection{Algorithmic Parameters and Basic Experimental Setup}
\label{sec:Setup}
In what follows, we use a custom implementation of MCTS written in C++ and use the mixed-integer optimization software Gurobi 5.0 \cite{gurobi} to solve the MO formulation.  All experiments were conducted on a computational grid with 2.2 GHz cores with 4GB of RAM in a single-threaded environment.  Although it is possible to parallelize many of the computations for each of the four methods, we do not explore this possibility in these experiments.  

Many of our experiments will study the effect of varying various hyperparameters (like the time limit per iteration) on the solution quality.  Unless otherwise specified in the specific experiment, all hyperparameters are set to their baseline values in Table~\ref{tab:BaselineParams}.
\begin{table}[!h]
\caption{Default parameters for algorithms \label{tab:BaselineParams}}
\begin{center}
\begin{tabular}{l l r r}
\toprule
Method   &   	Parameter   &   	Value   &   \\
\midrule
MCTS   &   	Time Limit per Iteration   &   	60\,s   &   
\\
   &   	Exploration Bonus $c$   &   	50   &   
\\
   &   	Rollout Policy   &   	FW Heuristic   &   
\\
    &   	Rollout Length   &   	10   &   
\\
    &   	Progressive Widening, Action Space $\alpha$   &   	0.5   &   
\\
    &   	Progressive Widening, State Space $\alpha^\prime$   &   	0.2   &   
\\
    &   	Progressive Widening, Action Space $k$   &   	40   &   
\\
    &   	Progressive Widening, State Space $k^\prime$   &   	40   &   

\\
    &   	Algorithm~\ref{alg:getnextga} $(u^\prime, u^{\prime \prime})$   &   	(0.3, 0.3)   &   
 \\
MO   &   	Time Limit per Iteration   &   	60\,s   &   
\\
    &   	Horizon Length   &   	10   &   
\\
\bottomrule
\end{tabular}
\end{center}
\end{table}

To ease comparison in what follows, we generally present the performance of each our algorithms relative to the performance of a randomized suppression heuristic.  At each time step, the randomized suppression heuristic chooses $| \mathcal{I} |$ cells without replacement from those cells which are currently burning and assigns suppression teams to them.  This heuristic should be seen as a naive ``straw man'' for comparisons only.  We will also often include the performance of our more tailored heuristic, the Floyd-Warshall (FW) heuristic, as a more sophisticated straw man. 

There are two experimental setups that we will return to throughout this section. 

\subsubsection{Grid 1}
In this setup, we consider a $k \times k$ grid with a varying reward function. There is a negative one reward received when the lower left cell is burning and the reward for a cell burning increases by one when traversing up or to the right across the grid. Also, the reward in the upper right hand corner is always $-10$. Figure~\ref{fig:grid1-reward} shows the rewards for a $k=8$ grid. 
\begin{figure}[!ht]
\centering
	\begin{tikzpicture}
	\begin{axis}[width=7cm,height=7cm,view={0}{90},nodes near coords,
	nodes near coords align={anchor=center},
	grid=major,
	xmin=1,xmax=9,ymin=1,ymax=9,
	xtick={1,2,...,9},
	ytick={1,2,...,9},
	xticklabels={},
	yticklabels={}]
		\addplot3[only marks] table {experiment1-reward-8.dat};
	\end{axis}
	\end{tikzpicture}
\caption{Rewards for Grid 1 with $k=8$.}\label{fig:grid1-reward}
\end{figure}

The fire in this experiment propagates as described in Section~\ref{sec:Dynamics} with 
\[
P(x,y) = \begin{cases} 0.06 & \text{ if } y \in \mathcal{N}(x)
\\
0 & \text{ otherwise.}
\end{cases}
\]
We also assume for this experiment that suppression efforts are successful with an 80\% probability---that is, $Q(x) =0.8$ for all $x \in \mathcal{X}$.   

For a single simulation we randomly generate an initial fire configuration---that is, whether or not each cell is burning and the fuel level in each cell. After generating an initial fire, suppression then begins according to one of our four approaches with $| \mathcal{I} |$ teams. The simulation and suppression efforts continue until the fire is extinguished or the entire area is burned out.   A typical experiment will repeat this simulation many times with different randomly generate initial fire configurations and aggregate the results. 

\begin{table}[!ht]
\caption{Grid 1 initial fire statistics}\label{tab:gridone-stats}
\begin{center}
\begin{tabular}{l r r r r r}
\toprule
 & $k=8$ & $k=12$ & $k=16$ & $k=20$ & $k=30$ \\
\midrule
Mean cells burning & 37.6 & 91.4 & 168.7 & 275.5 & 664.2 \\
Maximum cells burning & 62 & 142 & 244 & 372 & 845 \\
Mean fuel level for burning cells & 15.8 & 19.9 & 22.8 & 25.7 & 31.4 \\
Fuel level for non-burnt cells & 24 & 29 & 34 & 38 & 46 \\
\bottomrule
\end{tabular}
\end{center}
\end{table}

The specific process for initializing the fire configuration in a simulation is as follows:
\begin{enumerate}
\item Initialize all of the cells with a fuel level of $\lfloor \frac{k}{2P(x,y)} \rfloor$ and seed a fire in the lower left hand cell.
\item Allow the fire to randomly propagate for $\lfloor \frac{k}{2P(x,y)} \rfloor$ steps. Note that the lower left hand cell will naturally extinguish at the point.
\item Next, scale the fuel levels by a factor of $k^{-0.25}$. We scale the fuel levels to reduce the length of experiments where the number of suppression teams is insufficient to successfully fight the fire.
\end{enumerate}
Table~\ref{tab:gridone-stats} shows summary statistics for the initial fire configurations that arise from this process for this experiment.  We stress that the primary goal of this experiment is to explore the scalability of the various approaches.

\subsubsection{Grid 2}
This setup mirrors the setup in the previous experiment with two exceptions. First, we initialize fires in the middle of the grid.  Second, the reward function for cells is exponential across the grid.

\begin{figure}[!ht]
\centering
\includegraphics[width=.75\textwidth]{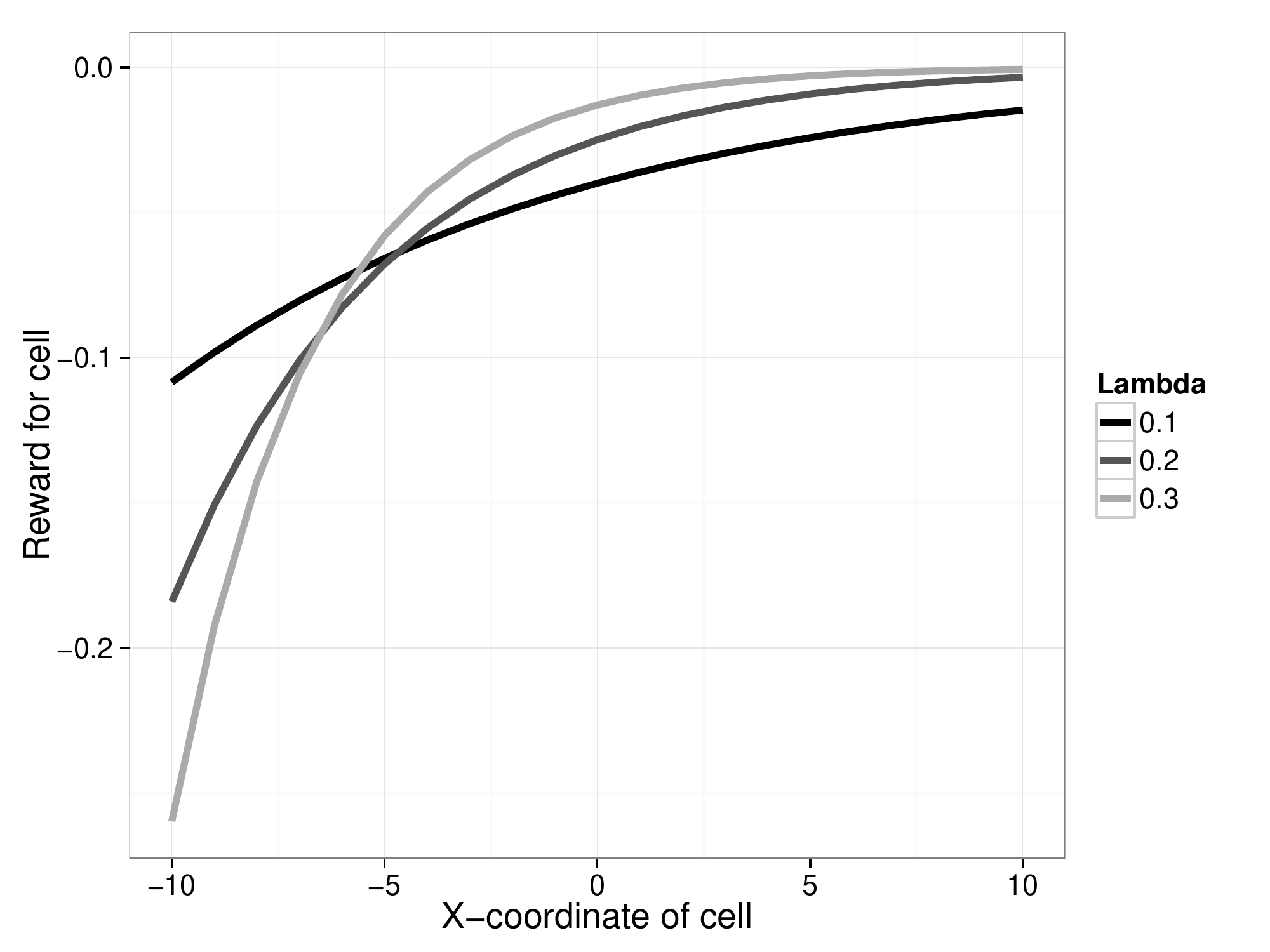}
\caption[Grid 2 reward structure.]{\label{fig:Grid2Lambda} Grid 2 reward structure for various values of $\lambda$.}
\end{figure}

Specifically, at time $t=0$ we ignite the cell in the middle of the grid, i.e., the cell at location $(\lceil k/2 \rceil, \lceil k/2 \rceil)$.  
The reward for cell $x = (i, j)$ is $-C \cdot \exp(-\lambda i)$ where $C^{-1} \equiv \sum_{i=1}^k e^{-\lambda i}$.  
Notice that the reward only depends on the horizontal location of the cell in the grid.  
In other words, cells located to the left are more valuable.  
The value of $\lambda$ controls the rate at which the reward grows.  
Some typical reward curves are shown in Figure~\ref{fig:Grid2Lambda} for a $k= 20$ grid.  
Observe that for large values of $\lambda$, the local reward structure at the site of the fire may seem quite flat.  
Good policies need to account for the fact that despite this local structure, suppression of cells on the righthand side of the fire is ultimately more valuable than suppression to the left.  

In this experiment, suppression efforts are still 80\% successful. The spread probabilities are given by
\[
P(x,y) = \begin{cases} 0.02 & \text{ if } y \in \mathcal{N}(x)
\\
0 & \text{ otherwise.}
\end{cases}
\]

As in the previous experiment, we begin with a random initial fire configuration.  
The specific process for initializing fire situations in this experiment is as follows:
\begin{enumerate}
\item Initialize the cells with fuel levels of $\lfloor \frac{k}{4P(x,y)} \rfloor$ and seed the fire in the cell in the middle of the grid---that is, the cell at location $(\lceil k/2 \rceil, \lceil k/2 \rceil)$.
\item Allow the fire to randomly spread for $\lfloor \frac{k}{4P(x,y)} \rfloor$ steps.  Note that we selected smaller fuel levels and initial times for the fire to build because the fire now starts in the middle of the grid.
\item Next, scale the fuel levels by a factor of $k^{-0.25}$. 
\end{enumerate}

Table~\ref{tab:gridtwo-stats} shows summary statistics for the initial fire configurations randomly generated in this experiment.
This experiment was designed to explore the ability of the various approaches to consider how the allocation of suppression teams now will impact future reward received. Because the reward function may be very different outside the local region of the fire, the approaches may need to plan many steps into the future.

\begin{table}[!ht]
\caption{Grid 2 initial fire statistics\label{tab:gridtwo-stats}}
\begin{center}
\begin{tabular}{l r r r}
\toprule
 & $k=9$ & $k=17$ & $k=25$ \\
\midrule
Mean cells burning & 37.8 & 154.9 & 363.7 \\
Maximum cells burning & 69 & 224 & 487 \\
Mean fuel level for burning cells & 5.3 & 7.5 & 8.9 \\
Fuel level for non-burnt cells & 8 & 11 & 13 \\
\bottomrule
\end{tabular}
\end{center}
\end{table}

\subsection{Tuning Hyperparameters for the MCTS Methodology}
\label{sec:HyperParamsMCTS}

\begin{figure}
\centering
\begin{subfigure}{\textwidth}
\includegraphics[width=\textwidth]{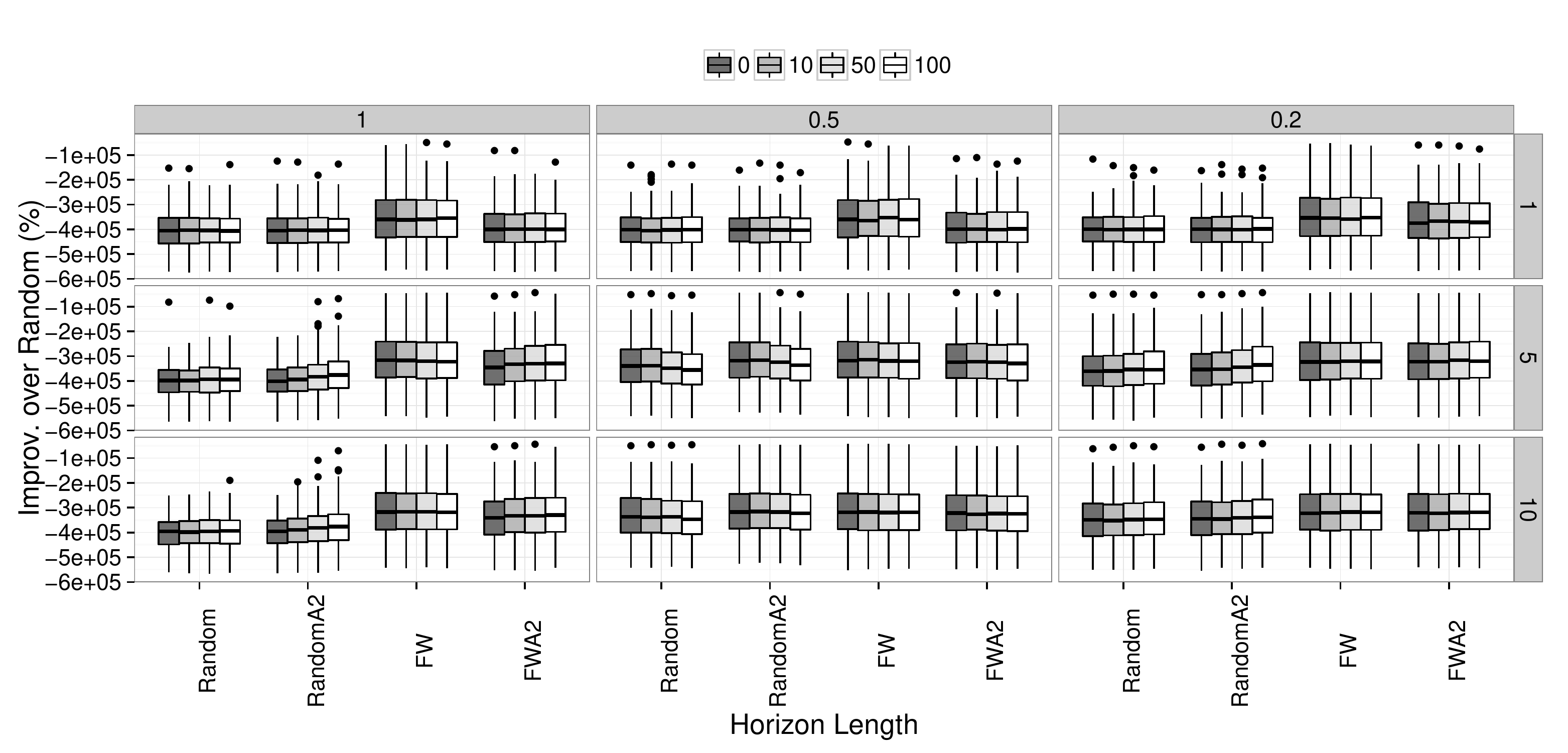}
\caption{}
\end{subfigure}
\\
\begin{subfigure}{\textwidth}
\includegraphics[width=\textwidth]{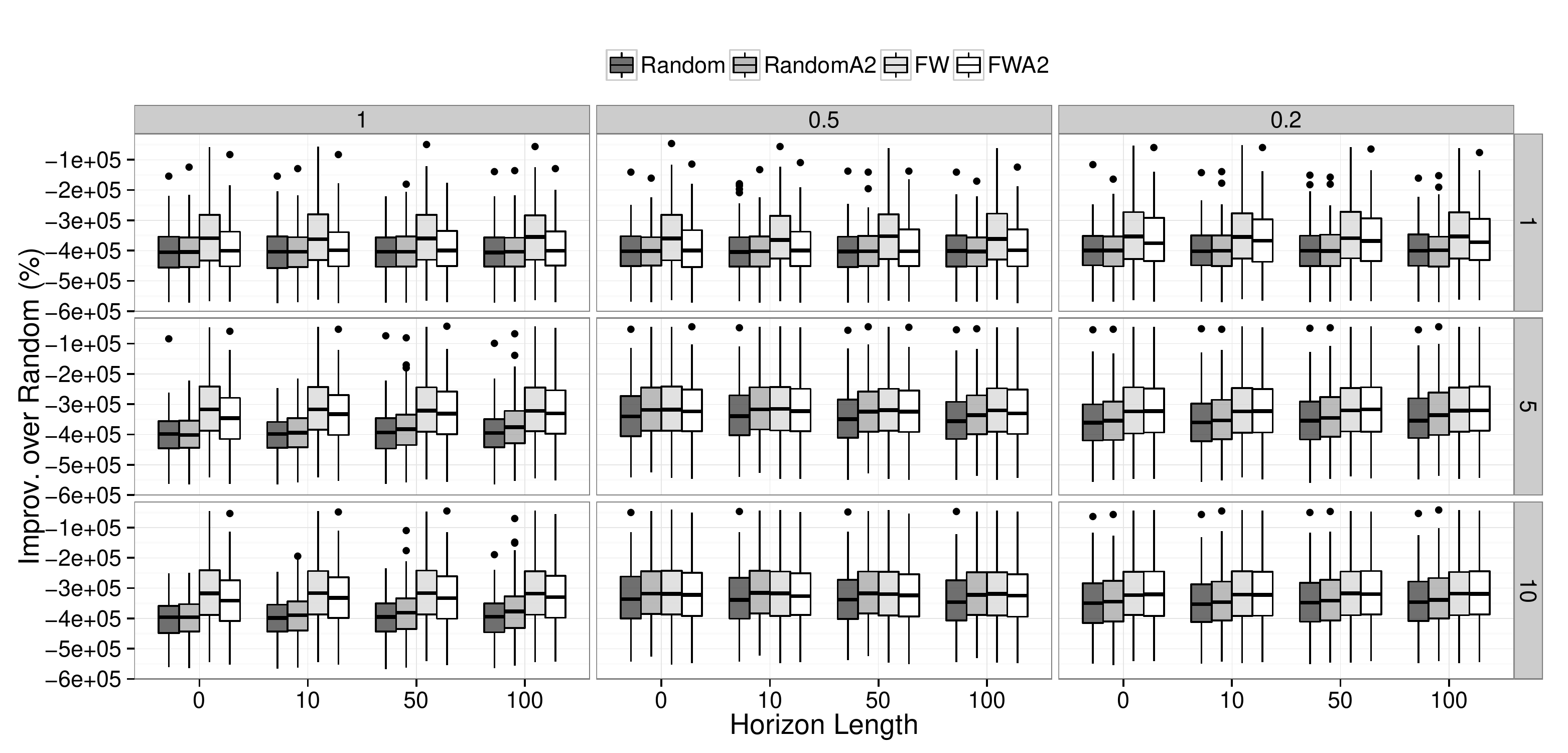}
\caption{}
\end{subfigure}
\caption[Hyperparameters of the MCTS algorithm.]{\label{fig:MCTSParams} The cumulative reward of the MCTS algorithm for various values organized by the exploration parameter $c$ and the rollout heuristic.  ``Random'' and ``FW'' refer to the random burn and Floyd-Warshall heuristics.  ``A2'' indicates that we additionally use our Algorithm~\ref{alg:getnextga} (based on genetic algorithm search heuristics) in the action generation phase.}
\end{figure}

The MCTS approach includes a number of hyperparameters to control the performance of the algorithm.  In this section, we explore the effects of some of these parameters using Grid 1, with $k=20$ and $4$ assets, and a time limit of 10\,s per iteration.  We vary the  exploration bonus $c \in \{ 0, 10, 50, 100\}$, the progressive widening factors $\alpha = \alpha^\prime \in \{1, 0.5, 0.2\}$, the depth of the rollout tree $d \in \{1, 5, 10\}$, whether or not we use Algorithm~\ref{alg:getnextga} (based on genetic algorithm search heuristics) in action generation, and whether we use the  random suppression heuristic or the FW heuristic in the rollout.  For each combination of hyperparameters, we run $256$ simulations and aggregate the results.
Figure~\ref{fig:MCTSParams} presents box plots of the cumulative reward. The top panel groups these box plots by the exploration bonus $c$, while the bottom one groups the same plots by the heuristic used.  Several features are noticeable:
\begin{enumerate}
	\item The effect of the exploration bonus $c$ is small and depends on the heuristic used.  For the FW heuristic, the effect is negligible.  One explanation for this features is that the FW heuristic is fairly strong on it own.  Hence, the benefits from local changes to this policy are somewhat limited. For the random suppression heuristic, there may be some effect, but it seems to vary with $\alpha$.  When $\alpha=0.2$, increased exploration benefits the random suppression heuristic, but when $\alpha=0.5$, less exploration is preferred.  
	\item From the second panel, in all cases it seems that MCTS with the FW heuristic outperforms MCTS with the random suppression heuristic.  The benefit of Algorithm~\ref{alg:getnextga} in action generation, however, is less clear.  Using Algorithm~\ref{alg:getnextga} does seem to improve the performance of the random suppression heuristic in some cases.  For the FW heuristic, there seems to be little benefit.  Again, one explanation of this phenomenon is that the FW heuristic already generates fairly good actions on its own.  With only 10\,s of computational time, it is difficult for the MCTS algorithm to find superior alternatives.
\end{enumerate}

To assess the statistical significance of these differences, we fit an additive effects model \citep{rice2007mathematical}.  To simplify the analysis, we fit two separate models, one for the random suppression heuristic and one for the FW heuristic.  Moreover, to simplify the models we use backward stepwise variable deletion beginning with the a full model with interactions of order two \citep{trevor2001elements}.  See Table~\ref{tab:MCTSParams} for the results.

\begin{table}
\centering
\caption{ \label{tab:MCTSParams} Estimated effects for the MCTS hyperparameters}
 \begin{tabular}{ l D{.}{.}{-1}D{.}{.}{-1} D{.}{.}{-1} D{.}{.}{-1}} 
\toprule 
  & \multicolumn{ 2 }{ c }{ FW Heuristic } & \multicolumn{ 2 }{ c }{ Random Suppression } \\ \midrule

 (Intercept)                    &    	 -353903.78 ^{***} &    	(0.00)&    	 -408969.13 ^{***}&    	 (0.00)              \\
$\alpha=0.5$                   &    	-3073.68&    	(0.37) &    	 5335.56          &    	 (0.16)              \\
$\alpha=0.2$                   &    	 7292.33 ^*        &    	(0.03) &    	 4616.83          &    	 (0.23)              \\
$\text{Depth}=5$                 &    	 34551.16 ^{***}   &    	(0.00) &    	 4211.96          &    	 (0.14)              \\
$\text{Depth}=10$                &    	 35375.98 ^{***}   &    	(0.00)&    	 3952.83          &    	 (0.17)              \\
A2                     &    	 -40434.84 ^{***}  &    	(0.00)&    	       -976.72            &    	    (0.71)                 \\
$c=10$                    &    	&	&    	 2857.04          &    	 (0.32)              \\
$c=50$                    &    	&	&    	 6900.73 ^*       &    	 (0.02)              \\
$c=100$                   &    	&	&    	 9366.90 ^{**}    &    	 (0.00)              \\
$\alpha=0.5$ : $\text{Depth}=5$    &    	4412.72&    	(0.29) &    	 58653.80 ^{***}  &    	 (0.00)              \\
$\alpha=0.2$ : $\text{Depth}=5$    &    	 -11279.75 ^{**}   &    	(0.01) &    	 41290.86 ^{***}  &    	 (0.00)              \\
$\alpha=0.5$ : $\text{Depth}=10$   &    	2989.4&    	(0.48) &    	 65456.71 ^{***}  &    	 (0.00)              \\
$\alpha=0.2$ : $\text{Depth}=10$   &    	 -11282.11 ^{**}   &    	(0.01) &    	 47508.01 ^{***}  &    	 (0.00)              \\
$\alpha=0.5$ : A2        &    	 8467.03 ^*        &    	(0.01) &    	        6960.33 ^*          &    	          (0.02)           \\
$\alpha=0.2$ : A2        &    	 20363.68 ^{***}   &    	(0.00)&    	     -1457.66             &    	        (0.61)             \\
$\text{depth}=5$ : A2      &    	 23627.58 ^{***}   &    	(0.00)&    	&	   \\
$\text{depth}=10$ : A2     &    	 24100.29 ^{***}   &    	(0.00)&    	&	   \\
$\alpha=0.5$ : $c=10$       &    	&	&    	 -2543.39         &    	 (0.53)              \\
$\alpha=0.2$ : $c=10$       &    	&	&    	 -983.12          &    	 (0.81)              \\
$\alpha=0.5$ : $c=50$       &    	&	&    	 -10139.50 ^*     &    	 (0.01)              \\
$\alpha=0.2$ : $c=50$       &    	&	&    	 -1250.75         &    	 (0.76)              \\
$\alpha=0.5$ : $c=100$      &    	&	&    	 -16684.02 ^{***} &    	 (0.00)              \\
$\alpha=0.2$ : $c=100$      &    	&	&    	 -674.51          &    	 (0.87)              \\
$\text{depth}=5$ : A2  &    	&	&    	 12733.89 ^{***}  &    	 (0.00)              \\
$\text{depth}=10$ : A2 &    	&	&    	 12608.70 ^{***}  &    	 (0.00)              \\
 $R^2$                          & 0.06       &       & 0.14      &       \\ 
adj. $R^2$                     & 0.06       &       & 0.14    &         \\ 
 \midrule
\multicolumn{3}{l}{\footnotesize{$^\dagger$ significant at $p<0.10$; $^* p<0.05$; $^{**} p<0.01$; $^{***} p<0.001$}} 
\\
\multicolumn{3}{p{9cm}}{\footnotesize{See Section~\ref{sec:HyperParamsMCTS} for details. Baseline should 
be interpreted as the value of $\alpha=1$, Depth of 1, $c=0$, without Algorithm~\ref{alg:getnextga} (based on genetic algorithm search heuristics).}}\\
\bottomrule
\end{tabular} 
 \end{table}

One can check that the features we observed graphically from the box plots are indeed significant.  Moreover, the random suppression heuristic demonstrates interesting second order effects between the depth and $\alpha$.  The performance improves substantially when the depth is greater than one and we restrict the size of the searched actions.  One explanation is that both parameter serve to increase the quality of the search tree, i.e., its depth, and the accuracy of the estimate at each of the searched nodes.

\subsection{State Space Size}
\label{sec:StateSpace}
We first study the performance of our algorithms as the size of the state space grows.  We simulate the performance of each of our methods on Grid 1 with either $4$ or $8$ suppression teams, using our default values of the hyperparameters and varying $k \in \{8, 12, 16, 20, 30\}$.  For each algorithm and combination of parameters, we simulate $256$ runs and amalgamate the results.  

\begin{figure}
\centering
\begin{subfigure}{0.4\textwidth}
	\includegraphics[width=\textwidth]{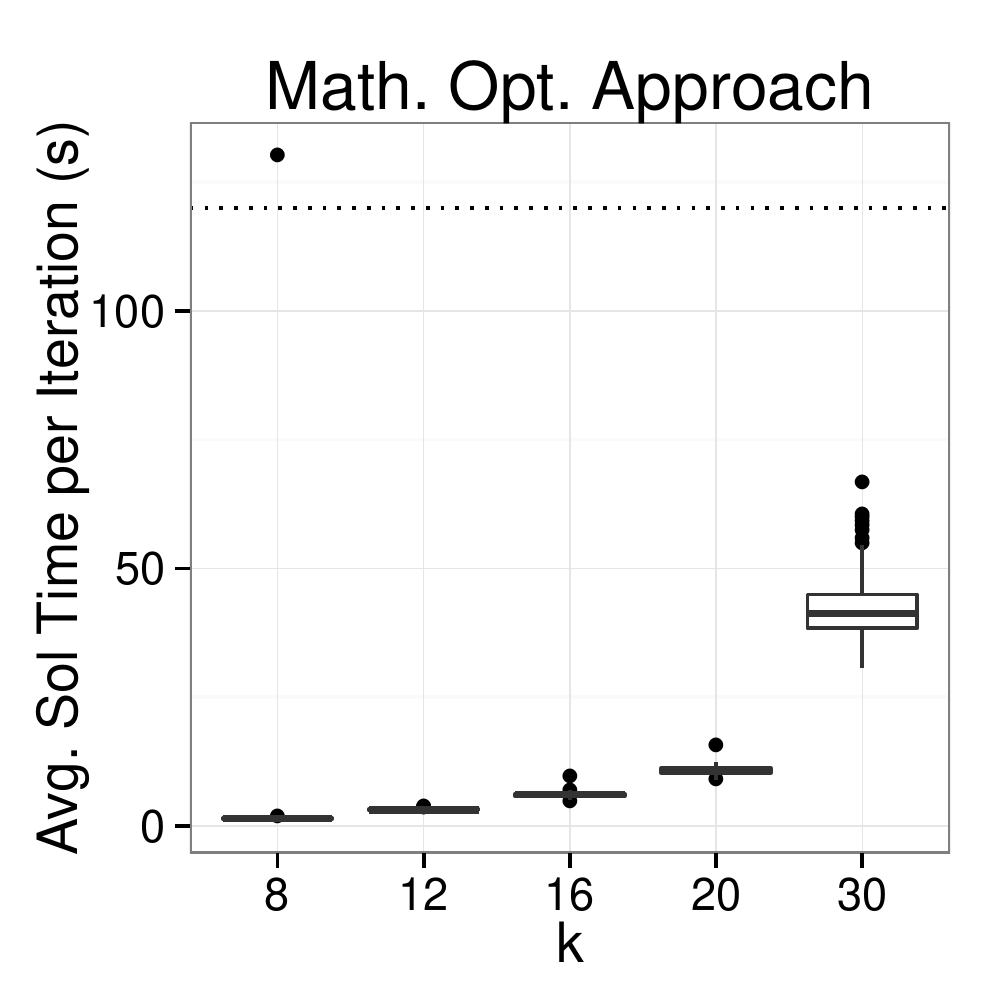}
	\caption{\label{fig:AvgSolTimes} Average}
\end{subfigure}
\begin{subfigure}{0.4\textwidth}
	\includegraphics[width=\textwidth]{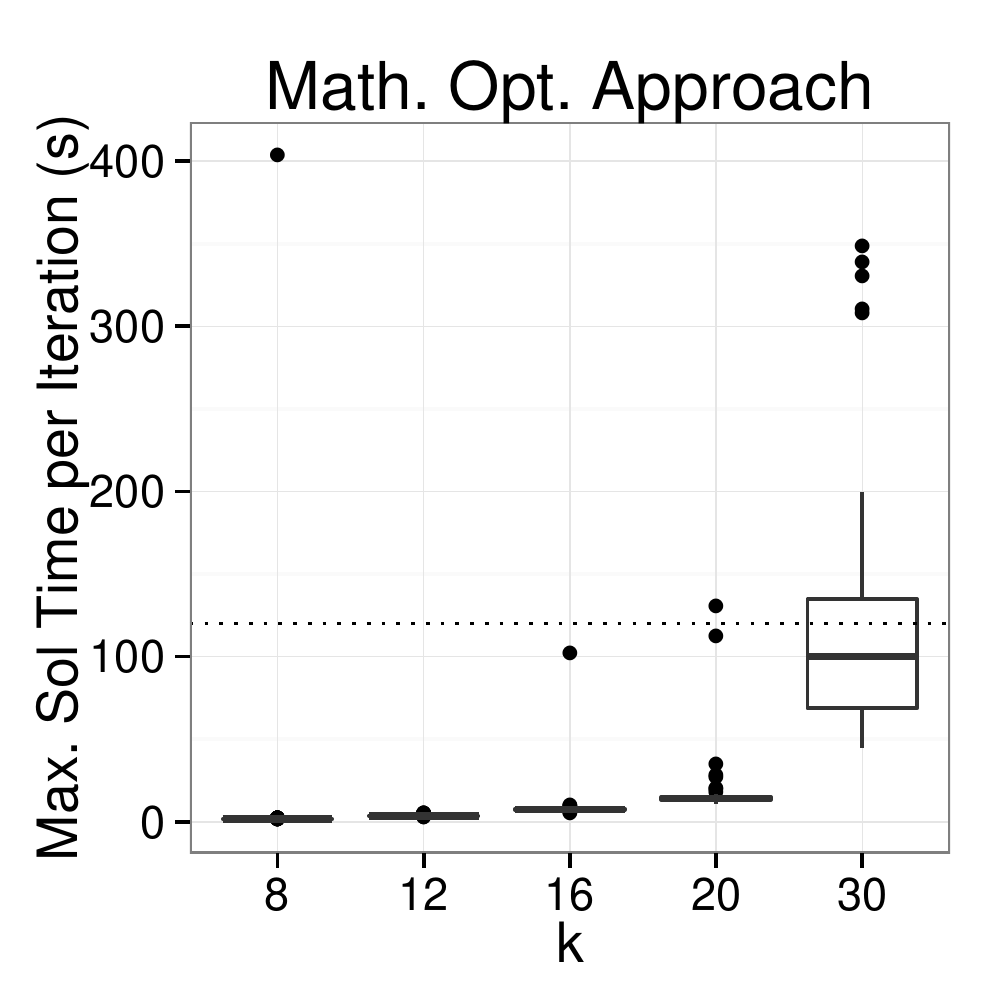}
	\caption {\label{fig:MaxSolTimes} Maximum}
\end{subfigure}
\caption[Solution time per iteration for the MO approach.] {Average and maximum iteration solution time with 8 teams and a desired time limit of 120\,s (dotted line).  Instances which exceed their allotted time are typically not solved to optimality.}
\end{figure}

Figures~\ref{fig:AvgSolTimes} and \ref{fig:MaxSolTimes} show the average and maximum solution time per iteration of the MO methodology when requesting at most 120\,s of computation time.  Notice that for most grids, the average time is well below the threshold -- in these instances the underlying integer program is solved to optimality.  For some grids, though, there are a few iterations which require much longer to find a feasible integer solution (cf.~the long upper-tail in Figure~\ref{fig:MaxSolTimes}).  Consequently, we compare our MO method to the MCTS method with both 60\,s and 120\,s of computation time.    

\begin{figure}
\centering
\begin{subfigure}{\textwidth}
	\includegraphics[width=\textwidth]{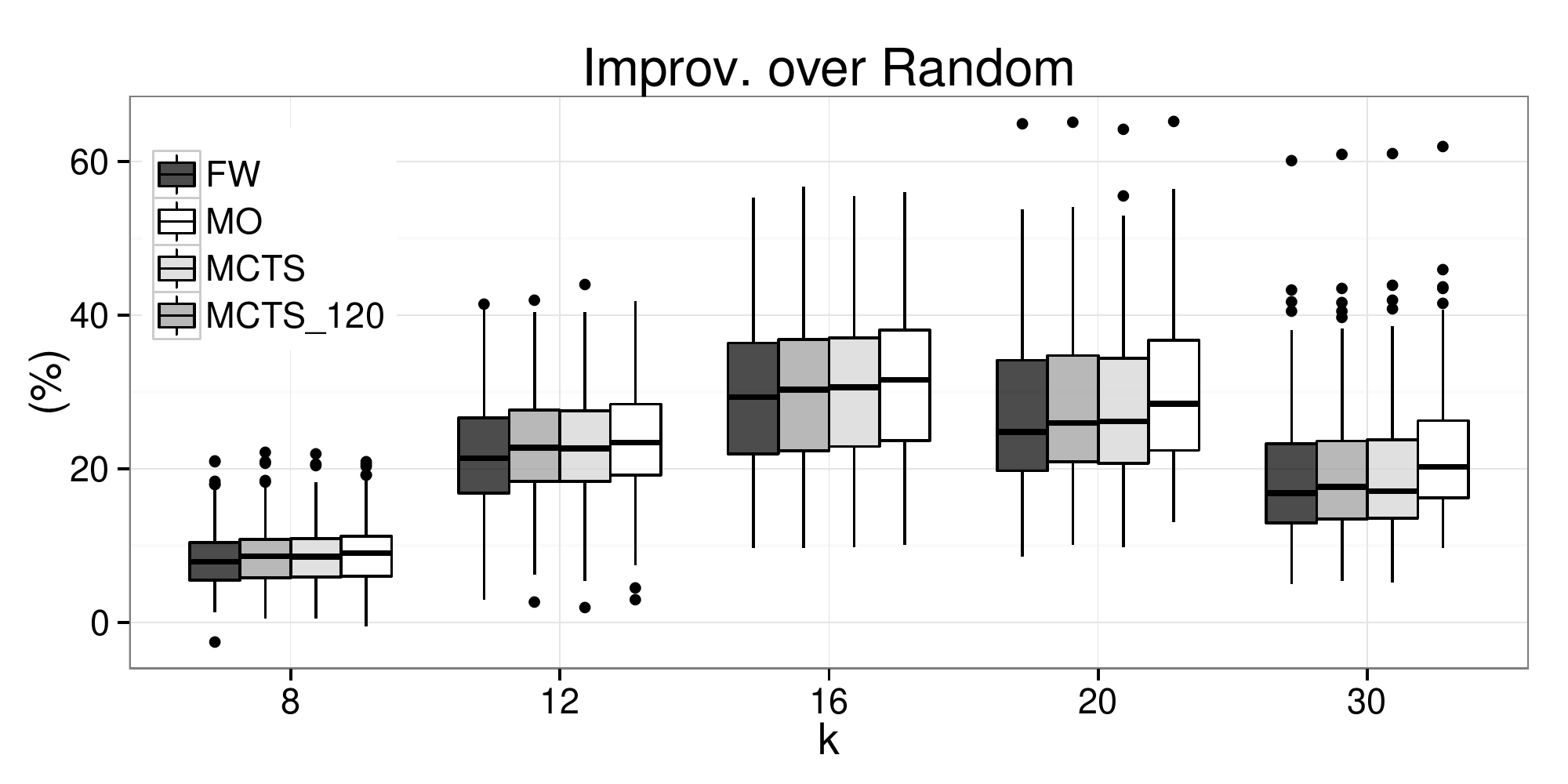}
	\caption{ \label{fig:Exp1Perfs8} 8 suppression teams.}
\end{subfigure}
\\
\begin{subfigure}{\textwidth}
	\includegraphics[width=\textwidth]{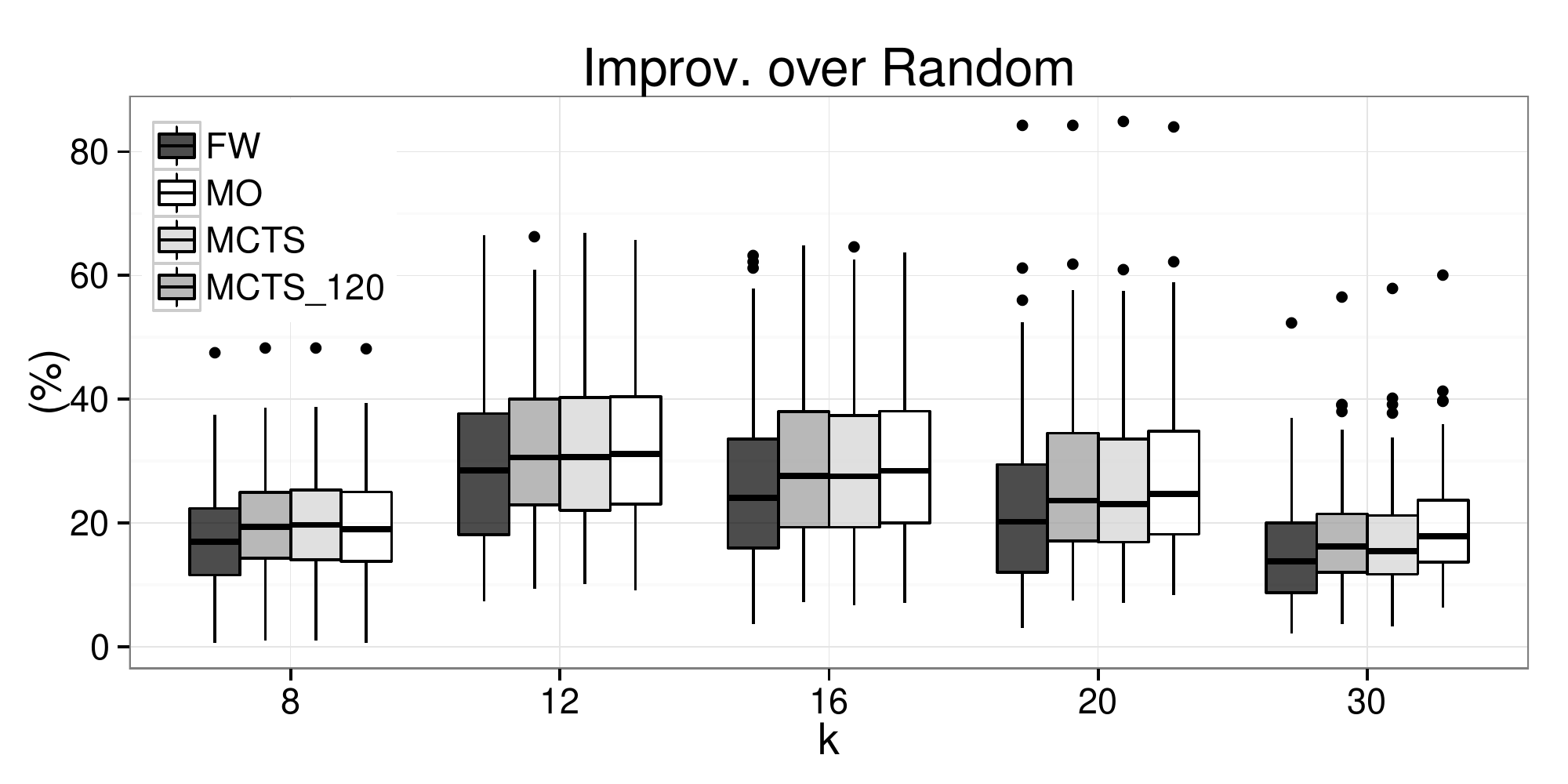}
	\caption{ \label{fig:Exp1Perfs4} 4 suppression teams.}
\end{subfigure}
\caption{Performance as a function of state space size.}
\end{figure}

A summary of the results is seen in Figures~\ref{fig:Exp1Perfs8} and \ref{fig:Exp1Perfs4}.  We stress that the runs with $4$ suppression teams are \emph{more} difficult than with $8$ teams; these instances are \emph{more} resource constrained.
Several features are evident in the plot.  First, all three methods seem to outperform the FW heuristic, but there seems to only be a small difference between the two MCTS runs.  The MO method does seem to outperform MCTS method, especially for tail-hard instances.  Namely, the lower whisker on the MO plot is often shorter than the corresponding whisker on the MCTS plots.  

To assess some of the statistical significance of these differences, we again fit two additive effects model; one for $8$ suppression teams and one for $4$.  In both cases, there are no significant second order interactions.    The coefficients of the first order interactions and corresponding $p$-values are shown in Table~\ref{tab:StateSpace}.  The values suggest that the two MCTS methods are very similar, with a slight edge for the 120\,s variant, and that the MO method with a time limit of 60\,s does outperform both.  
\begin{table}
\centering
\caption{\label{tab:StateSpace} Estimated effects for the percentage improvement relative to the random heuristic }
\begin{tabular}{ l D{.}{.}{-1}D{.}{.}{-1} D{.}{.}{-1} D{.}{.}{-1} } 
\toprule
  & \multicolumn{ 2 }{ c }{ 8 Teams } & \multicolumn{ 2 }{ c }{ 4 Teams } \\ 
  & \multicolumn{ 1 }{ c }{ Coefficient } & \multicolumn{ 1 }{ c }{ $p$-value } &   \multicolumn{ 1 }{ c }{ Coefficient } & \multicolumn{ 1 }{ c }{ $p$-value } \\ 
\midrule
(Intercept)           & 7.77 ^{***}  & (0.00) & 16.92 ^{***} 
                             & (0.00)      \\ 
$k=12$           & 14.46 ^{***} & (0.00) & 11.98 ^{***} 
                             & (0.00)      \\ 
$k=16$           & 21.89 ^{***} & (0.00) & 9.75 ^{***}  
                             & (0.00)      \\ 
$k=20$           & 19.53 ^{***} & (0.00) & 6.09 ^{***}  
                             & (0.00)      \\ 
$k=30$           & 10.99 ^{***} & (0.00) & -2.12 ^{***}  
                             & (0.00)      \\ 
MCTS (120\,s) & 0.87 ^*      & (0.03) & 3.03 ^{***}  
                             & (0.00)      \\ 
MCTS (60\,s)  & 0.83 ^*      & (0.04) & 2.74 ^{***} 
                             & (0.00)      \\ 
MO  & 2.22 ^{***}  & (0.00) & 3.80 ^{***}  
                             & (0.00)       \\
 $R^2$                 & 0.47         & & 0.21        \\ 
adj. $R^2$            & 0.47         & & 0.20        \\ 
\midrule
\multicolumn{3}{l}{\footnotesize{$^\dagger$ significant at $p<0.10$; $^* p<0.05$; $^{**} p<0.01$; $^{***} p<0.001$}} \\
\multicolumn{5}{p{9cm}}{\footnotesize{For details, see Section~\ref{sec:StateSpace}.  The intercept should
be interpreted as baseline of $k=8$ with the FW heuristic.}}\\
\bottomrule
\end{tabular} 
 \end{table}

\begin{figure}
\centering
	\includegraphics[width=\textwidth]{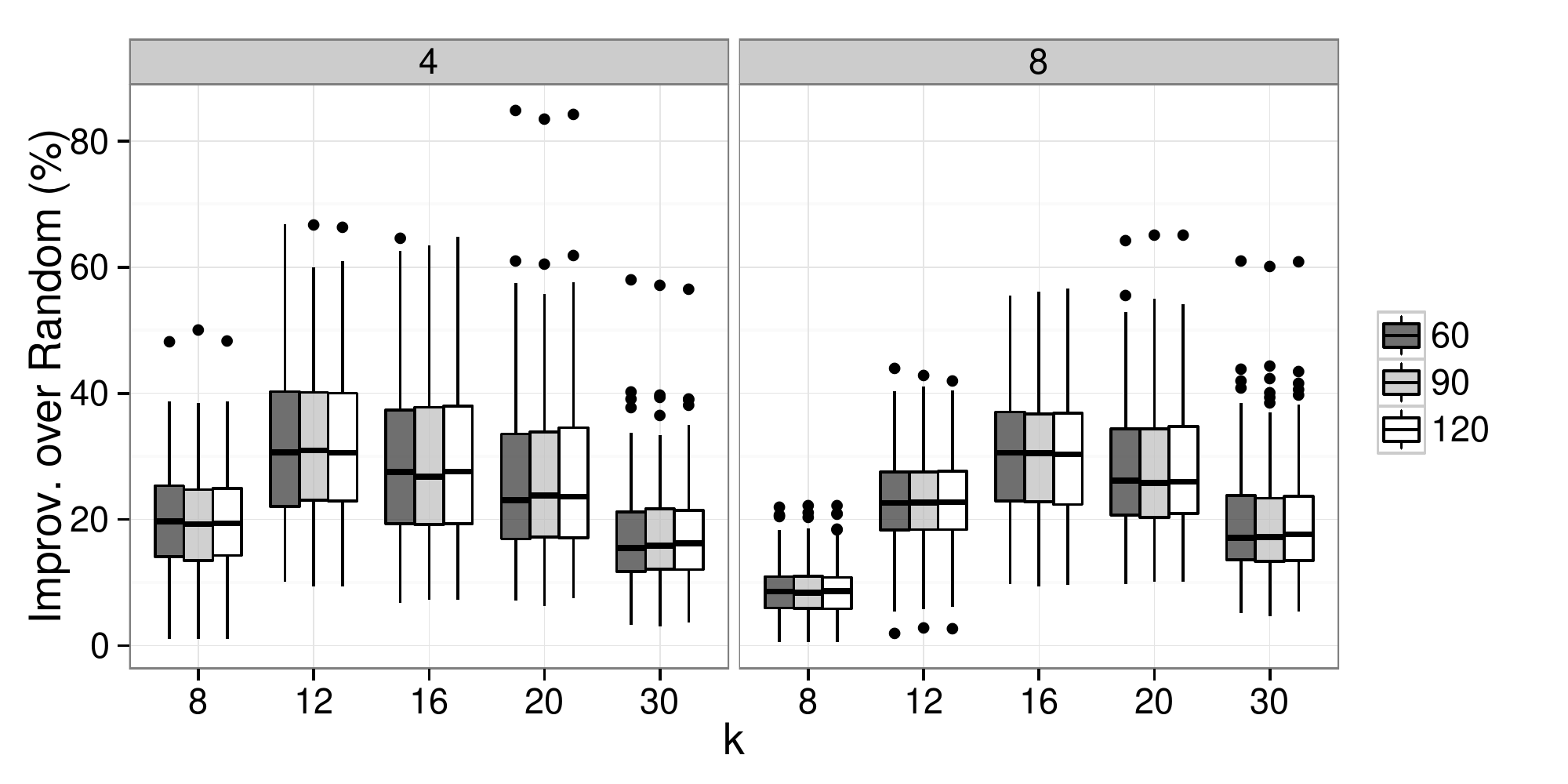}
	\caption[MCTS performance as a function of time allotted.]{ \label{fig:MCTSByTime}  MCTS performance as a function of time allotted.  See also Section~\ref{sec:StateSpace}.}
\end{figure}

To further understand the effect on performance of varying the time limit per iteration of the MCTS method, we re-ran the above experiment with 60\,s, 90\,s and 120\,s of time for the MCTS method.  Results are summarized in Figure~\ref{fig:MCTSByTime}.
As can be seen, the difference in performance is small.  

In summary then, these results suggest that MCTS and MO perform comparably for small state spaces, but as the size of the state space grows, MO begins to outperform MCTS, with the magnitude of the edge growing with the state space size.


\subsection{Action Branching Factor}
\label{sec:ActionBranch}
Intuition suggests that the performance of the MCTS algorithm is highly dependent on the magnitude of 
the action branching factor, i.e., the number of actions available from any given state.  As mentioned in the main text, without progressive widening, when the action branching factor is larger than the number of iterations, the MCTS algorithm will only expand the search tree to depth one.  Even with progressive widening, choosing good candidate actions is critical to growing the search tree in relevant directions.  
A simple calculation using Stirling's approximation confirms that for the MDP outlined in Section~\ref{sec:Dynamics}, the action branching factor at time $t$ is given by 
\begin{equation}
 \frac{N_B(t)^{| \mathcal{I}|} }{{| \mathcal{I}|}!} \approx \frac{\left(e \cdot  \frac{N_B(t)}{| \mathcal{I}|} \right)^{| \mathcal{I}|}}{\sqrt{2 \pi {| \mathcal{I}|}}},  \label{eq:stirling}
\end{equation}
where $N_B(t)$ is the number of cells that are burning at time $t$.  For even medium sized-instances, this number is extremely large.
Consequently, in this section we study the performance of our algorithms with respect to the action branching factor.

We have already seen initial results in Section~\ref{sec:HyperParamsMCTS} suggesting that both our progressive widening and Algorithm~\ref{alg:getnextga} for action generation improve upon the base MCTS algorithm.  It remains to see how MCTS with these refinements compares to our mathematical optimization formulation.  We compute the relative improvement of the MCTS and MO approaches over the randomized suppression heuristic on Grid 1 over $256$ simulations and $k=10$.  
\begin{figure}
\centering
\includegraphics[width=\textwidth]{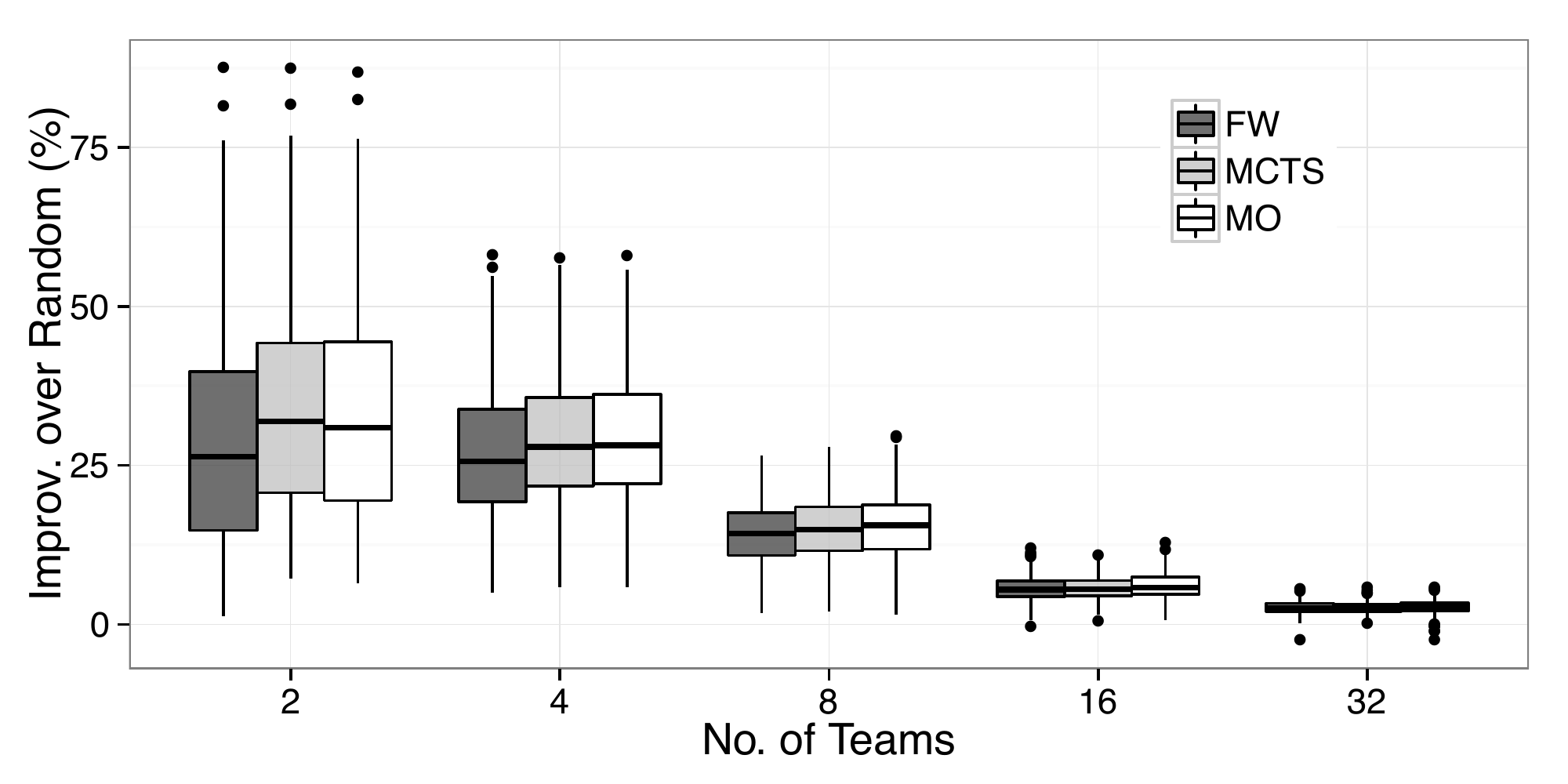}
\caption{Performance as function of number of suppression teams, $k=$ 10. \label{fig:PerfByTeamsSmall}}
\end{figure}
Figure~\ref{fig:PerfByTeamsSmall} summarizes the average relative improvement for each our methods.

\begin{figure}
\centering
\begin{subfigure}{.4\textwidth}
\includegraphics[width=\textwidth]{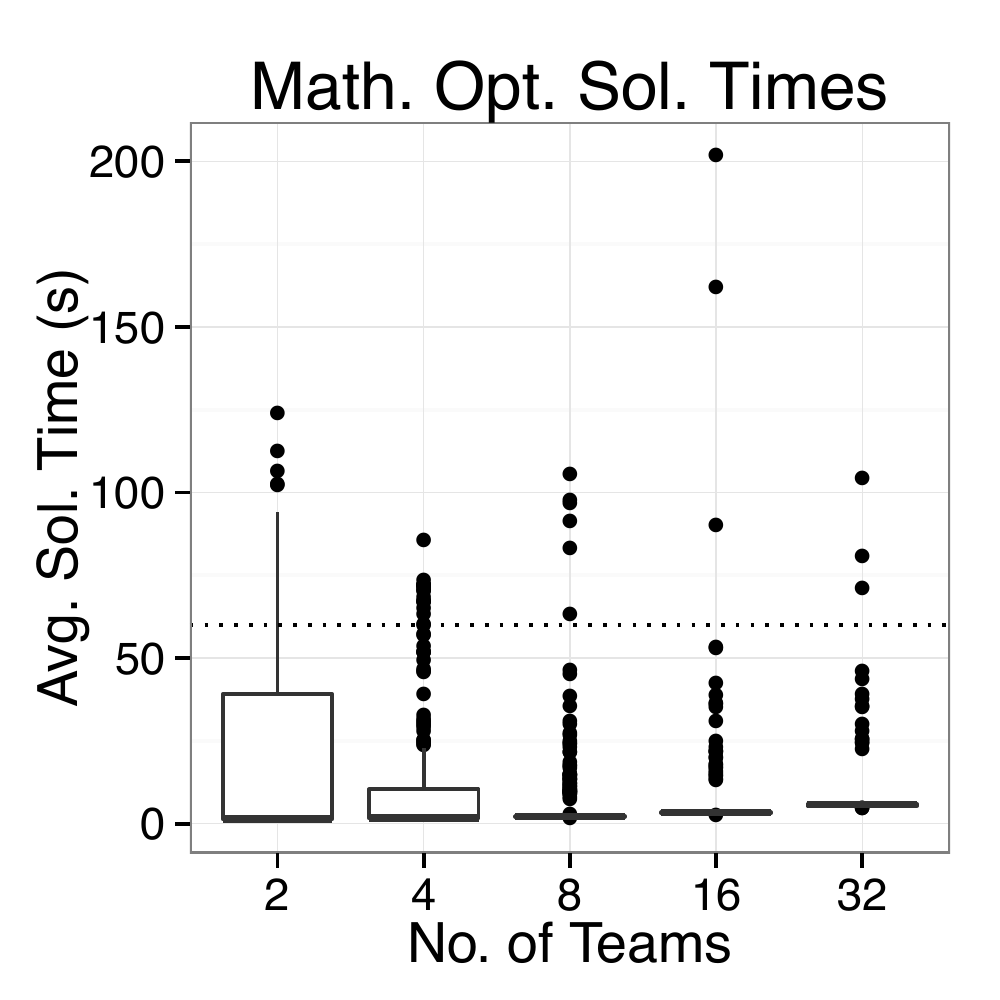}
\caption{ \label{fig:SolTimesbyTeamSmall} $k=10$}
\end{subfigure}
\begin{subfigure}{.4\textwidth}
\includegraphics[width=\textwidth]{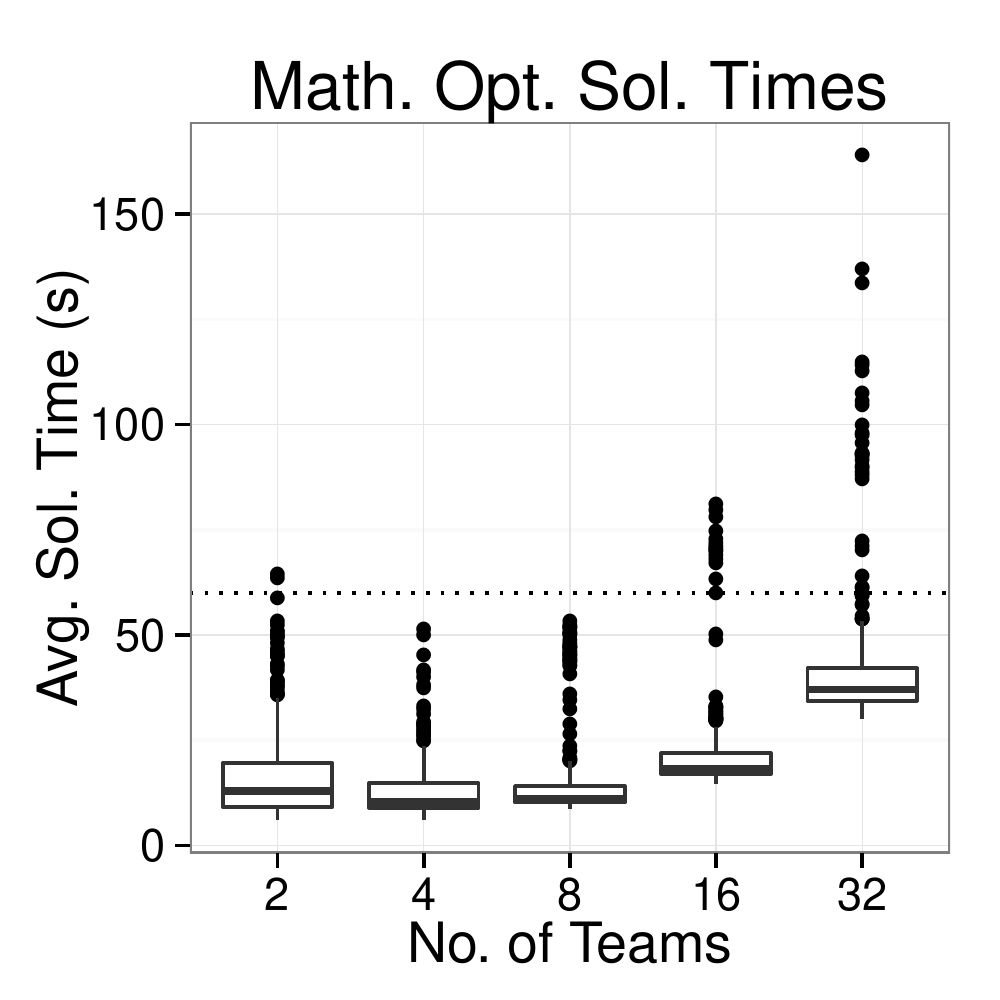}
\caption[Solution times per iteration by number of suppression teams, MO approach]{ \label{fig:SolTimesbyTeam} $k=20$}
\end{subfigure}
\caption{MO average solution times.}
\end{figure}

Recall that it is not possible to control for the exact time used by the MO algorithm. Figure~\ref{fig:SolTimesbyTeamSmall} shows a box-plot of the average time per iteration for the MO approach.  Based on this plot, we feel that comparing the results to the MCTS algorithm with 60\,s of computational time is a fair comparison.  

\begin{figure}
\centering
\includegraphics[width=\textwidth]{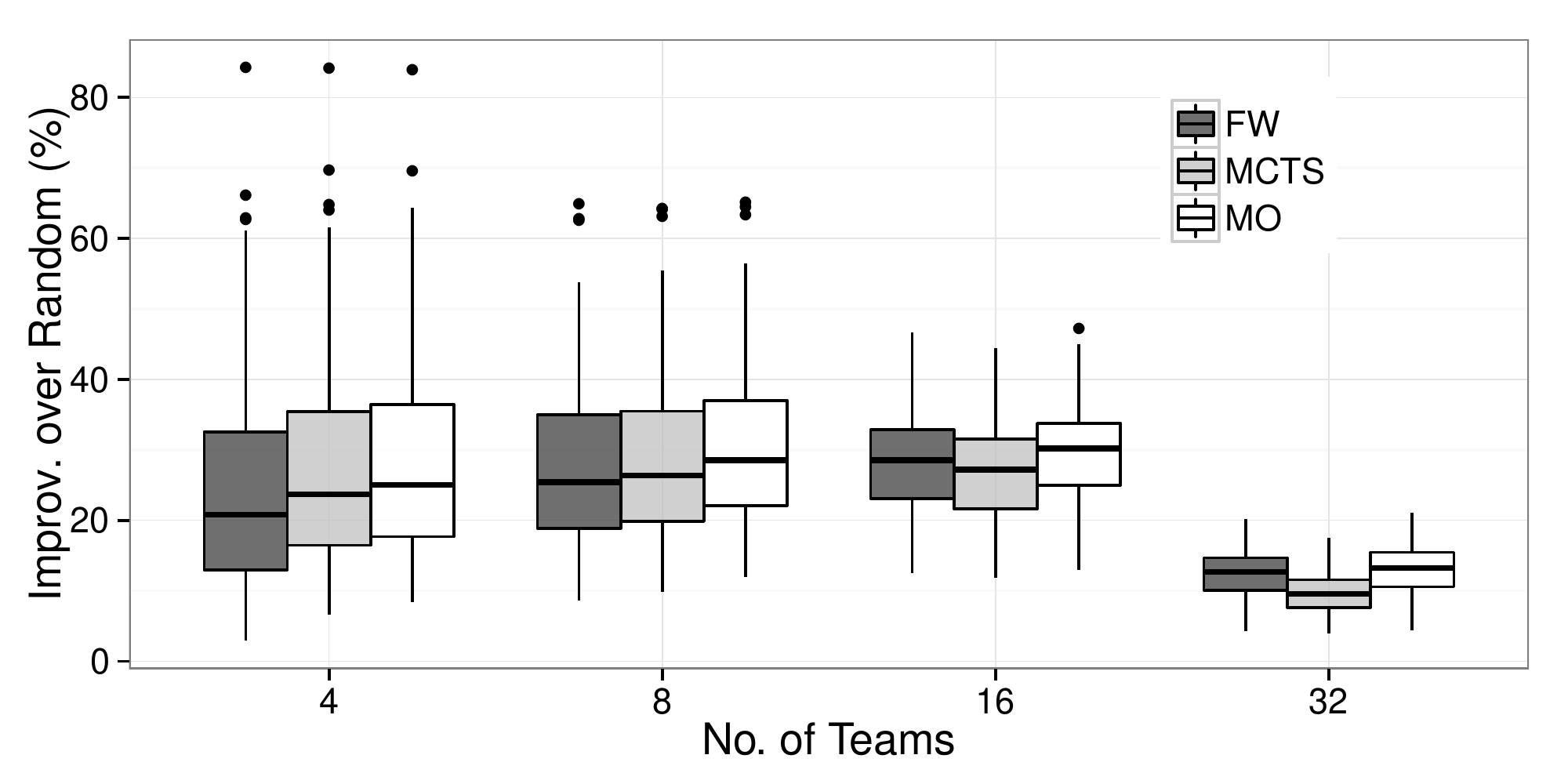}
\caption{Performance as a function of number of suppression teams, $k=$ 20. \label{fig:PerfByTeams}}
\end{figure}

The relative performance of all our methods degrades as the number of teams become large.  This is principally because the randomized suppression heuristic improves with more teams. Although the FW heuristic is clearly inferior, the remaining approaches appear to perform similarly.  Indeed, ANOVA testing suggests the differences between these to MO and MCTS, the difference between MO and MCTS is not statistically significant ($p$-values are well above $0.4$).  To try and isolate more significant differences between the methodologies, we re-run the above experiment with $k=20$.  The results can be seen in Figure~\ref{fig:PerfByTeams} and the average solution times in Figure~\ref{fig:SolTimesbyTeam}.  

In contrast with the previous experiment, MO appears to outperform MCTS.  Interestingly, although MCTS seems to outperform the FW heuristic for a small number of suppression teams, it performs worse for more teams.  To test the significance of these differences, we fit a linear regression model for the improvement over the randomized suppression heuristic as a function of the number of teams, the algorithm used, and potential interactions between the number of teams and the algorithm used.  The results are in Table~\ref{tab:AdditiveEffects}.  The intercept value is the baseline, \emph{fitted value of the FW heuristic policy for 4 teams.}  Several observations can be drawn from Table~\ref{tab:AdditiveEffects}.  First, MO outperforms FW for all choices of team size with statistical significance.  MCTS, however, is statistically worse than FW with 16 or 32 teams.

\begin{table}
\begin{center}
	\caption[Estimated effects of team size]{ \label{tab:AdditiveEffects} Estimated effects for MCTS and MO with varying team size }
	\begin{tabular}{ l D{.}{.}{-1} D{.}{.}{-1}} 
	\toprule
	 & \multicolumn{ 1 }{ c }{ Coefficients} &  \multicolumn{ 1 }{ c }{ $p$-value} \\
	\midrule
	(Intercept)                            & 23.40 ^{***} 
	                                       & (0.00)          \\ 
	8 Teams              & 3.92 ^{***}    
	                                       & (0.00)          \\ 
	16 Teams             & 4.79 ^{***}      
	                                       & (0.00)          \\ 
	32 Teams             & -10.84 ^{***}   
	                                       & (0.00)          \\ 
	MCTS                            & 3.01 ^{***}      
	                                       & (0.00)          \\ 
	MO                            & 4.55 ^{***}      
	                                       & (0.00)          \\ 
	8 Teams : MCTS  & -1.93 ^\dagger  
	                                       & (0.10)          \\ 
	16 Teams : MCTS & -4.36 ^{***}     
	                                       & (0.00)          \\ 
	32 Teams : MCTS & -5.86 ^{***}     
	                                       & (0.00)          \\ 
	8 Teams : MO & -1.49           
	                                       & (0.20)          \\ 
	16 Teams : MO  & -3.19 ^{**}         
	                                       & (0.01)          \\ 
	32 Teams : MO  & -4.02 ^{***}    
	                                       & (0.00)           \\
	$R^2$                                  & 0.36            \\ 
	adj. $R^2$                             & 0.36            \\ 
	\midrule
	\multicolumn{2}{l}{\footnotesize{$^\dagger$ significant at $p<0.10$; $^* p<0.05$; $^{**} p<0.01$; $^{***} p<0.001$}} \\
	\multicolumn{3}{p{9cm}}{\footnotesize{See also Section~\ref{sec:ActionBranch}.  The intercept should be
	interpreted as the performance of four teams under the FW algorithm.}}
	\\
	\bottomrule
\end{tabular} 
\end{center}
\end{table}

In summary, differences between the MO and MCTS methods become visible only when the grid size $k$ is large, i.e., when the instances are sufficiently ``difficult'' to solve.  It appears that although progressive widening and Algorithm~\ref{alg:getnextga} for action selection partially address the challenges of a large action state branching factor, the mathematical optimization approach is better suited to these instances.

\subsection{Asymmetric Costs and Horizon Length}
\label{sec:Exp2}
In this section we examine the performance of our algorithms on Grid 2 (cf.~Section~\ref{sec:Setup}).  Recall, in this setup, the cost structure is asymmetric; cells to the right side of the grid are more valuable than cells to the left.  At the same time, the local reward structure at the point of ignition is relatively flat.  A good algorithm must recognize the differential value despite the local reward structure.  

\begin{figure}
\centering
\includegraphics[width=\textwidth]{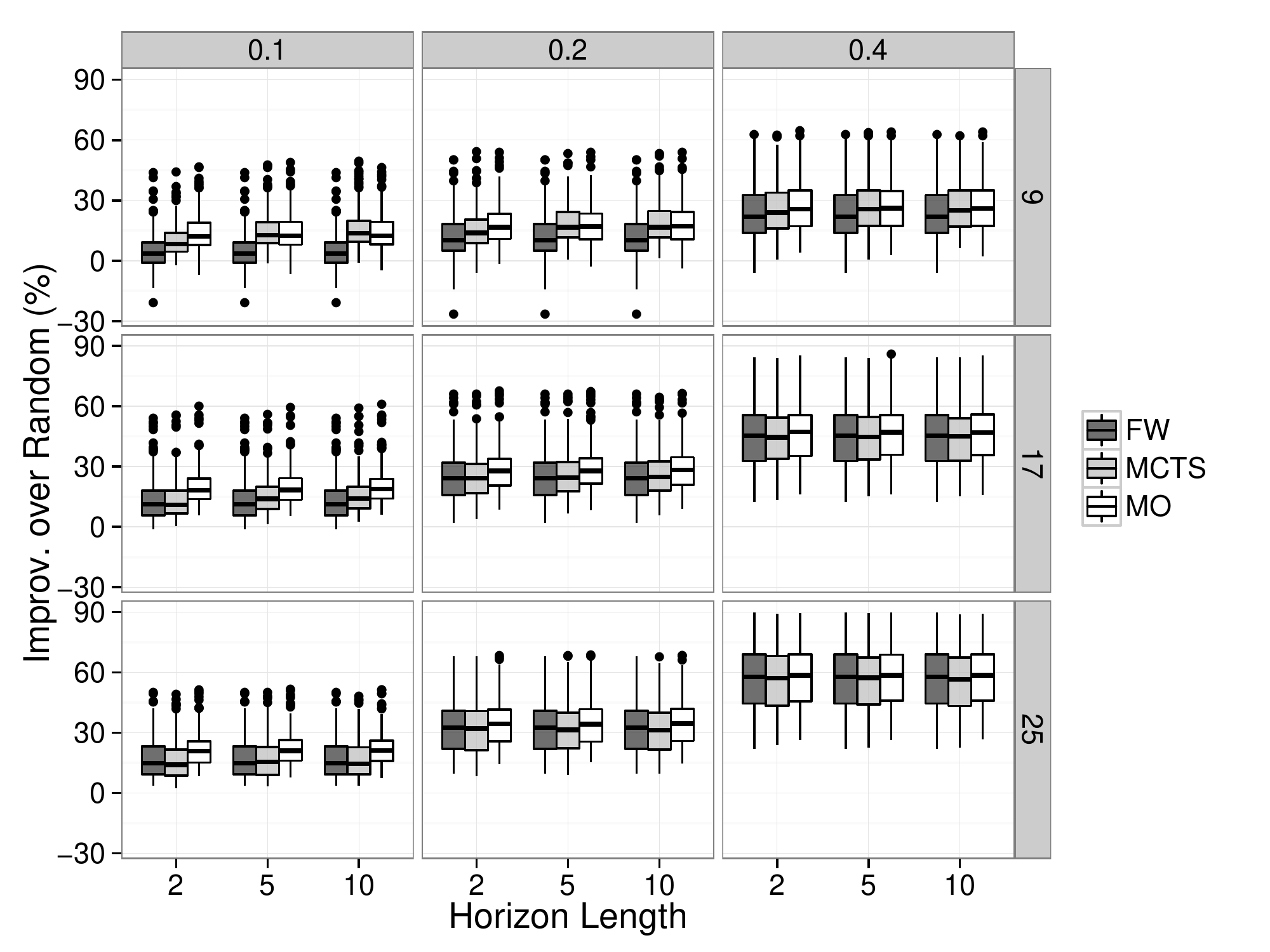}
\caption[Performance as a function of horizon length.]{ \label{fig:PerfBylambda} Performance as a function of horizon length.  Different panels correspond to different values of $k \in \{$9, 17, 25$\}$ and and $\lambda \in \{$0.1, 0.2, 0.4$\}$.  }
\end{figure}

\begin{figure}
\centering
\includegraphics[width=\textwidth]{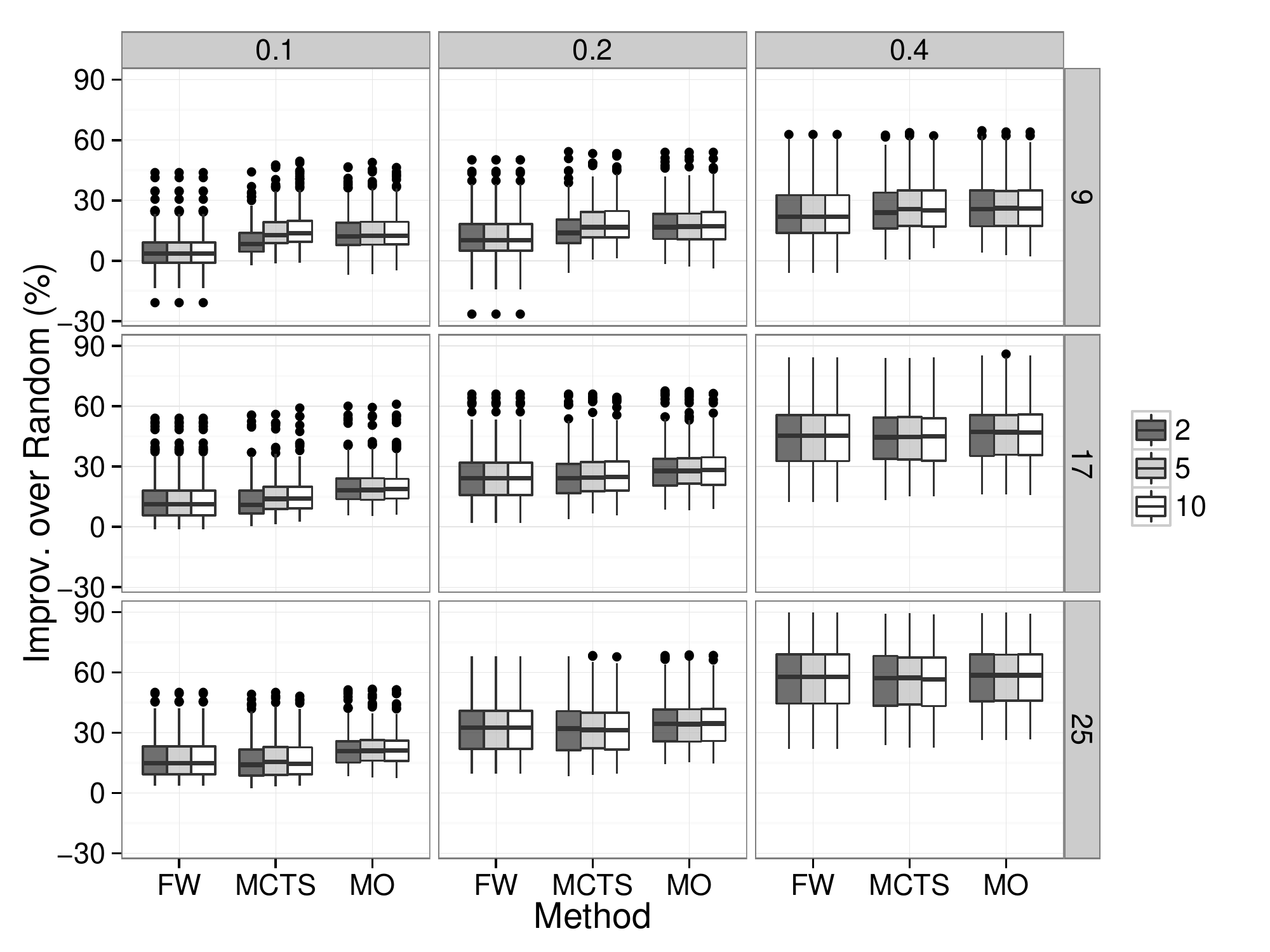}
\caption[Performance as a function of horizon length by method.]{ \label{fig:PerfBylambdaMethod} Performance as a function of horizon length, grouped by method.  Different panels correspond to different values of $k \in \{$9, 17, 25$\}$ and and $\lambda \in \{$0.1, 0.2, 0.4$\}$.  }
\end{figure}

We consider different combinations of $\lambda \in \{0.1, 0.2, 0.4\}$, $k \in \{9, 17, 25\}$, and horizon length in $\{2, 5, 10\}$.  For each combination, we consider $256$ simulations.  Figure~\ref{fig:PerfBylambda} summarizes the performances of each of our methods.  For all the methods, there is an upward trend as $\lambda$ increases.  For small values of $\lambda$ (i.e., flatter reward structure), MO has a marked edge over the other methods.  As $\lambda$ increases, the difference shrinks.  
Figure~\ref{fig:PerfBylambdaMethod} shows the same box plots, but grouped by method.  This plot suggests that increasing the horizon length from $2$ to $5$ improves performance, but there is negligible improvement for length $10$. This agrees to some extent with our intuition: as the horizon length increases, the algorithms become less myopic and are better able to account for the rapid growth in reward as one moves to the right of the grid.

We investigate the statistical significance of these differences by fitting additive effects models.  The full model contains many insignificant interactions.  We drop insignificant variables using backwards stepwise deletion.  The resulting fit can be seen in Table~\ref{tab:Exp2Fitted}.  The fitted model suggests that differences observed in the previous plot are statistically significant.

\begin{table}
\begin{center}
\caption[Estimated effects for Section~\ref{sec:Exp2}] {\label{tab:Exp2Fitted} Estimated effects for Section~\ref{sec:Exp2}}
\begin{tabular}{ l D{.}{.}{-1} D{.}{.}{-1} } \toprule
Effect & \multicolumn{ 1 }{ c }{ Coefficient }  & \multicolumn{ 1 }{ c }{ $p$-value }\\ \midrule
(Intercept)                   & 5.45^{***} &   (0.00)  \\
$\lambda = 0.2$             & 6.66^{***}   & (0.00)   \\
$\lambda = 0.4$             & 17.70^{***}   & (0.00)  \\
$k = 17$                   & 7.84^{***}   & (0.00)  \\
$k = 25$                   & 10.87^{***}   & (0.00) \\
Depth 5                & 0.00       & (1.00)  \\
Depth 10               & 0.00        & (1.00)  \\
MCTS                     & 6.25^{***}   & (0.00) \\
MO                       & 8.76^{***}    & (0.00) \\
$\lambda = 0.2$ : $k = 17$ & 5.22^{***}   & (0.00)  \\
$\lambda = 0.4$ : $k = 17$ & 14.38^{***}   & (0.00) \\
$\lambda = 0.2$ : $k = 25$ & 9.33^{***}    & (0.00) \\
$\lambda = 0.4$ : $k = 25$ & 23.56^{***}   & (0.00) \\
$\lambda = 0.2$: MCTS   & -1.47^{**}   & (0.00) \\
$\lambda = 0.4$ : MCTS   & -2.94^{***}  & (0.00) \\
$\lambda = 0.2$ : MO     & -2.81^{***}    & (0.00) \\
$\lambda = 0.4$ : MO     & -4.98^{***}  & (0.00) \\
$k = 17$ : MCTS         & -5.14^{***}  & (0.00) \\
$k = 25$ : MCTS         & -6.51^{***}  & (0.00) \\
$k = 17$ : MO           & -2.33^{***}  & (0.00) \\
$k = 25$ : MO           & -3.48^{***}  & (0.00) \\
Depth 5 : MCTS      & 1.50^{**}   & (0.00) \\
Depth 10 : MCTS     & 1.47^{**}    & (0.00)  \\
Depth 5 : MO        & 0.18        & (0.72)  \\
Depth 10 : MO       & 0.14        & (0.78)  \\
\midrule
 \multicolumn{2}{l}{\footnotesize{$^\dagger$ significant at $p<0.10$; $^* p<0.05$; $^{**} p<0.01$; $^{***} p<0.001$}} \\
 \multicolumn{3}{p{9cm}}{\footnotesize{The intercept should be interpreted as 
the value of the FW heuristic when $\lambda=.1$, $k=9$ and the horizon length is 2.}}
 \\
 \bottomrule
\end{tabular}
\label{table:coefficients}
\end{center}
\end{table}

\section{Conclusion}
\label{sec:Conclusion}
In this study, we consider a dynamic resource allocation problem motivated by tactical wildfire management, in which fire spreads stochastically on a finite grid and the decision maker must allocate resources at discrete epochs to suppress it. We propose two different solution approaches: one based on Monte Carlo tree search, and one based on mathematical optimization.


Our study makes three broad methodological contributions. The first of these contributions is to the understanding of MCTS: to the best of our knowledge, this study is the first application of Monte Carlo tree search to high-dimensional dynamic resource allocation motivated by a real-world application, and our results suggest that MCTS holds promise in other large scale stochastic control applications. Our numerical results uncover some interesting insights into how MCTS behaves in relation to parameters such as the exploration bonus, the progressive widening parameters and others, as well as larger components such as the method of action generation and the rollout heuristic. Our results show that these components are highly interdependent and cannot be calibrated separately of each other---for example, our results show that the choices of action generation method and progressive widening factor become very important when the rollout heuristic is not strong on its own (e.g., the random suppression heuristic) but are less valuable when the rollout heuristic is strong to begin with (e.g., the Floyd-Warshall heuristic). These insights will be valuable for practitioners interested in applying MCTS to other problems.

The second broad methodological contribution of our study is to the understanding of mathematical optimization in dynamic resource allocation. The formulation that we consider is a continuous, deterministic one that avoids the combinatorial nature of the true fire dynamics, which are discrete and stochastic. Our results indicate that such approximations, in spite of how much they depart from the true dynamics, can lead to effective heuristics for large scale stochastic control problems. This insight is both important and not obvious.

The third broad methodological contribution is towards the understanding of the relative merits of MCTS and mathematical optimization. Our results show that while both methodologies exhibit comparable performance for smaller instances, for larger instances, the mathematical optimization approach exhibits a significant edge. Initial evidence suggests this edge may be related more closely to action branching factor than the state space branching factor. 



\section*{Acknowledgements}
This paper is the result of research and development sponsored by the Assistant Secretary of Defense for Research and Engineering, ASD(R\&E). The work of the fifth author was supported by a PGS-D award from the Natural Sciences and Engineering Research Council (NSERC) of Canada.

\bibliographystyle{elsarticle-harv}
\bibliography{references} 

\end{document}